\documentclass[article,hidelinks,onefignum,onetabnum]{siamart220329}

\usepackage{lipsum}
\usepackage{amsfonts}
\usepackage{graphicx}
\usepackage{epstopdf}
\usepackage{algorithmic}
\ifpdf
  \DeclareGraphicsExtensions{.eps,.pdf,.png,.jpg}
\else
  \DeclareGraphicsExtensions{.eps}
\fi

\newsiamremark{remark}{Remark}
\newsiamremark{hypothesis}{Hypothesis}
\crefname{hypothesis}{Hypothesis}{Hypotheses}
\newsiamthm{claim}{Claim}

\headers{Error bounds for rank-one DNN reformulations of QAP}{Y. T. Qian, S. H. Pan, S. J. Bi and H. D. Qi}

\title{Error Bounds for Rank-one Double Nonnegative Reformulations of QAP and Exact Penalties \thanks{Submitted to the editors DATE.
\funding{The second author's work was funded by the National Natural Science Foundation of China under project 12371299. The third author's work was funded by the National Natural Science Foundation of China under project 12371323. 
Qian and Qi were funded by Hong Kong RGC General Research Fund PolyU/15309223 and PolyU AMA Project 230413007.}}}

\author{Yitian Qian\thanks{Department of Data Science and Artificial Intelligence, The Hong Kong Polytechnic University, Hong Kong  
  (\email{yitian.qian@polyu.edu.hk}).}
\and Shaohua Pan\thanks{School of Mathematics, South China University of Technology, Guangzhou
(\email{shhpan@scut.edu.cn},
\email{bishj@scut.edu.cn}).}
\and Shujun Bi\footnotemark[3]
  \and Houduo Qi\thanks{Department of Data Science and Artificial Intelligence, and Department of Applied Mathematics, The Hong Kong Polytechnic University, Hong Kong 
  (\email{houduo.qi@polyu.edu.hk})}}

\usepackage{amsopn}

\ifpdf
\hypersetup{
  pdftitle={Error Bounds for Rank-one Double Nonnegative Reformulations of QAP and Exact Penalties},
  pdfauthor={Y. T. Qian, S. H. Pan, S. J. Bi and H. D. Qi}
}
\fi

\usepackage{url}
\usepackage{amsmath,amsfonts,amssymb,bm,color,caption}
\usepackage{cases}
\usepackage{multirow}
\usepackage{longtable}
\usepackage{makecell}
\usepackage{booktabs}
\usepackage{subfigure}
\usepackage{epstopdf}
\usepackage{rotating}
\usepackage{nccmath}
\usepackage{stmaryrd}

\begin{document}

\maketitle

\begin{abstract}
This paper focuses on the error bounds for several equivalent rank-one doubly nonnegative (DNN) conic reformulations of the quadratic assignment problem (QAP), a class of challenging combinatorial optimization problems. We provide three equivalent rank-one DNN reformulations of the QAP, including the one proposed in \cite{Jiang21}, and establish the locally and globally Lipschitzian error bounds for their feasible sets. Then, these error bounds are employed to prove that the penalty problems induced by the difference-of-convexity (DC) reformulation of the rank-one constraint are global exact penalties, and so are the penalty problems for their Burer-Monteiro (BM) factorizations. As a byproduct, the penalty problem for the rank-one DNN reformulation in \cite{Jiang21} is shown to be a global exact penalty without the calmness assumption. Finally, we illustrate the application of these exact penalties by proposing a relaxation approach with one of them to seek a rank-one approximate feasible solution. This relaxation approach is validated to be superior to the commercial solver Gurobi for \textbf{132} benchmark instances in terms of the relative gap between the generated objective value and the known best one and the number of instances with better objective values. 
\end{abstract}



\section{Introduction}\label{sec1.0}

 Quadratic assignment problem (QAP) is a fundamental one in location theory, which allocates $n$ facilities to $n$ locations and minimizes a quadratic function on the distance between the locations and the flow between the facilities. Let $\mathbb{R}^{n\times q}$ represent the space consisting of all $n\times q$ real matrices, equipped with the trace inner product $\langle\cdot,\cdot\rangle$ and its induced Frobenius norm $\|\cdot\|_F$. The QAP is formulated as
 \begin{equation}\label{QAP}
 \min_{X\in\mathcal{P}}\,\langle X,AXB+C\rangle,
 \end{equation}
 where $\mathcal{P}:=\{X\in\{0,1\}^{n\times n}\ |\ Xe=e,X^{\top}e=e\}$ is the set of all $n\times n$ permutation matrices, and $A,B\in\mathbb{R}^{n\times n}$ and $C\in\mathbb{R}^{n\times n}$ are the data matrices. It is worth pointing out that the set of global (or local) optimal solutions of \eqref{QAP} keeps unchanged when $A$ is replaced by the nonnegative matrix $A-\min_{i,j}A_{ij}ee^{\top}$ (similarly for $B$ and $C$).

 The QAP has wide applications in facility layout, chip design, scheduling, manufacturing, and so on; see \cite{Burkard13,Drezner15}. Notice that $\mathcal{P}=\{X\in\mathbb{R}_{+}^{n\times n}\,|\, X^{\top}X=I\}$, the set of nonnegative orthogonal matrices in $\mathbb{R}^{n\times n}$. When $C=0$, it can be rewritten as
 \begin{equation}\label{EQAP1}
 \min_{X\in\mathcal{P}}\|A^{\top}X+XB\|_F^2,
 \end{equation}
 which has a wide application in pattern recognition and machine vision \cite{Conte14}, and also covers the bandwidth minimization problem (see \cite[section 7]{JiangLW16}) appearing frequently in sparse matrix computations, circuit design, and VLSI layout \cite{Chinn82,Lai99}.

 Due to the curse of dimensionality, it is impractical to achieve a global optimal solution of the QAP; for example, when $n=30$, seeking a global optimal solution with the branch-and-bound method is still computationally challenging. This work aims at establishing locally Lipschitzian error bounds for the rank-one DNN conic reformulations of \eqref{QAP}, and using these error bounds to achieve their global exact penalties, one of which is applied to develop a continuous relaxation approach to yield an approximate feasible solution, so an upper bound for the optimal value.
 \subsection{Related works}\label{sec1.1}
 
 Consider that the forthcoming global exact penalties can be used to design relaxation methods for seeking an approximate feasible solution. We mainly review the relaxation methods to yield a feasible solution or an infeasible one whose objective value offers a lower bound for the optimal value. The first class is based on the problem obtained by replacing the discrete $\mathcal{P}$ with its convex hull $\mathcal{D}$: 
 \begin{equation}\label{relaxation1}
 \min_{X\in\mathcal{D}}\,\langle X,AXB+C\rangle.
 \end{equation}
 For example, Xia \cite{Xia10} proposed a Lagrangian smoothing algorithm by solving a sequence of $L_2$-regularization subproblems of \eqref{relaxation1} with dynamically reducing the regularization parameter. Observing that $\mathcal{P}=\mathcal{D}\cap\{X\!\in\mathbb{R}^{n\times n}\ |\ \|X\|_0=n\}$, where $\|X\|_0$ is the zero-norm (the number of nonzero entries) of $X$, Jiang et al. \cite{JiangLW16} considered the $L_0$-norm regularization problem of \eqref{relaxation1} that was proved to have the same global optimal solution set as \eqref{relaxation1}, and developed a regularization algorithm by replacing the zero-norm $\|X\|_0$ with its $p$-norm approximation for $p\in(0,1)$. Zaslavskiy et al. \cite{Zaslavskiy09} and Liu et al. \cite{LiuQiao12} proposed path-following algorithms by solving a sequence of convex combinations of convex and concave optimization problems over the set $\mathcal{D}$.

 The second one is developed by the following equivalent reformulation of \eqref{QAP}
 \begin{equation}\label{North}
 \min_{X\in\mathbb{R}^{n\times n}}\Big\{\langle A, (X\circ X)B(X\circ X)^{\top}\rangle+\langle C,X\rangle\ \ {\rm s.t.}\ \ X^{\top}X=I,\,X\in\mathbb{R}_{+}^{n\times n}\Big\},
 \end{equation} 
 where ``$\circ$'' denotes the Hadamard product. Wen and Yin \cite{Wen13} proposed a continuous relaxation method by applying the augmented Lagrangian method (ALM) to deal with the nonnegative constraint of \eqref{North} and solving every ALM subproblem with a nonmonotone line-search Riemannian gradient descent method. Qian et al. \cite{QianPanXiao22} proposed a continuous relaxation method by the global exact penalty of \eqref{North} induced by the Moreau envelope of the $\ell_1$-norm distance from the cone $\mathbb{R}_{+}^{n\times n}$, which solves a sequence of penalty problems with a nonmonotone line-search Riemannian gradient descent method. Different from the first one, this class of relaxation methods captures a certain surface information of $\mathcal{P}$ by leveraging the orthogonal manifold.

 The third one is proposed by the outer approximation to the completely positive cone $\mathcal{C}^*\!:=\!{\rm conv}\{zz^{\top}\,|\,z\in\mathbb{R}^{p}_+\}$ with $p=n^2$. Inspired by the work \cite{Anstreicher2000}, Povh and Rendl \cite{Povh09} proved that \eqref{QAP} is equivalent to the following completely positive conic program
 \begin{align}\label{Equad-CPP}
 &\min_{Y\!\in\mathbb{S}^p}\ \langle B\otimes A + {\rm Diag}({\rm vec}(C)),Y\rangle  \nonumber\\
 &\ {\rm s.t.}\ \ \langle I,Y^{ij}\rangle=\varpi_{ij}\ \ \forall i,j\in[n],\ {\textstyle\sum_{t=1}^n}{Y^{tt}}=I,\\
 &\qquad\ \langle ee^{\top},Y\rangle=p,\,Y\!\in\mathcal{C}^*,\nonumber
 \end{align}
 where ``$\otimes$'' denotes the Kronecker product, $\varpi_{ij}=1$ if $i=j$, otherwise $\varpi_{ij}=0$, and $Y^{ij}\in\mathbb{R}^{n\times n}$ for $i,j\in[n]$ denotes the $(i,j)$-th block of $Y\in\mathbb{S}^p$. Considering that the closed convex cone $\mathcal{C}^*$ is numerically intractable, they replaced $\mathcal{C}^*$ with the DNN matrix cone and obtained the following DNN conic convex relaxation of \eqref{Equad-CPP}
 \begin{align}\label{Rquad-DNN}
  &\min_{Y\in\mathbb{S}^{p}}\ \langle B\otimes A + {\rm Diag}({\rm vec}(C)),Y\rangle  \nonumber\\
  &\ {\rm s.t.}\ \ \langle I,Y^{ij}\rangle=\varpi_{ij}\ \ \forall i,j\in[n],\ {\textstyle\sum_{t=1}^n}Y^{tt}=I,\\
  &\qquad\ \langle ee^{\top},Y\rangle=p,\, Y\in K:=\mathbb{S}_{+}^p\cap\mathbb{R}_{+}^{p\times p}.\nonumber
 \end{align}
 Later, Kim et al. \cite{Kim16} considered a Lagrangian-DNN relaxation for another reformulation of \eqref{QAP} to  produce better lower bounds than those yielded by solving the DNN conic program \eqref{Rquad-DNN}. Recently, Oliveira et al. \cite{Oliveira18} and Graham et al. \cite{Graham22} applied the facial reduction technique to a DNN relaxation of \eqref{QAP}, a little different from \eqref{Rquad-DNN}, and developed a splitting method to solve the DNN relaxation. The above DNN relaxations benefit from the surface information of $\mathcal{P}$ by its lifted reformultion, and provide tighter lower bounds than positive semidefinite relaxations of \eqref{QAP}.

 The last one is developed by the rank-one DNN reformulation of \eqref{QAP}. Notice that the constraint sets of \eqref{Equad-CPP} and \eqref{Rquad-DNN} are different from that of \eqref{QAP} only in the rank-one constraint $Y\in\mathcal{R}:=\big\{Y\in\mathbb{S}^p\ |\ {\rm rank}(Y)\le 1\big\}$, which is equivalent to the DC constraint $\|Y\|_*-\|Y\|=0$. Hence, the problem \eqref{QAP} can be reformulated as 
 \begin{align}\label{EroneDNNcon}
 \min_{Y\in\mathbb{S}^{p}}\Big\{\langle G,Y\rangle\ \ {\rm s.t.}\ \ Y\in \Delta_0,\,\|Y\|_*-\|Y\|=0\Big\}
 \end{align}
 with $G:=\frac{1}{2}[(B\otimes A)+(B\otimes A)^{\top}]+{\rm Diag}({\rm vec}(C))$, where $\Delta_0$ is the set defined by 
 \begin{equation*}
  \Delta_0:=\Big\{Y\in K\,|\,\langle I,Y^{ij}\rangle=\varpi_{ij}\ \ \forall i,j\in[n],\ {\textstyle\sum_{t=1}^n}{Y^{tt}}=I,\langle ee^{\top},Y\rangle=p\Big\}.
 \end{equation*}
 The equivalence between \eqref{EroneDNNcon} and \eqref{QAP} reveals that the DC constraint becomes the hurdle for the solving of \eqref{QAP}. Since the handling of DC constraints is much harder than that of DC functions numerically, it is natural to focus on the penalty problem
 \begin{align}\label{EroneDNNcon-penalty}
  \min_{Y\in\mathbb{S}^{p}}\Big\{\langle G,Y\rangle+\rho\big(\|Y\|_*-\|Y\|\big)\ \ {\rm s.t.}\ \ Y\in\Delta_0\Big\},
 \end{align}
 and to prove that it is a global exact penalty of \eqref{EroneDNNcon}, i.e., there exists a threshold $\overline{\rho}>0$ such that the problem \eqref{EroneDNNcon-penalty} associated to every $\rho\ge\overline{\rho}$ has the same global optimal solution set as \eqref{EroneDNNcon}. The latter is significant for numerical computation of \eqref{EroneDNNcon}, since it means that computing a finite number of penalty problems \eqref{EroneDNNcon-penalty} with appropriately large $\rho$ can yield a favorable solution. Jiang et al. \cite{Jiang21} proved that the problem \eqref{EroneDNNcon-penalty} is a global exact penalty of \eqref{EroneDNNcon} by assuming that the multifunction
 \begin{equation}\label{Upsilon0-map}
  \Upsilon_0(\tau)\!:=\Big\{Y\in\Delta_0\ |\ \|Y\|_*-\|Y\|=\tau\Big\}\quad{\rm for}\ \tau\in\mathbb{R}
 \end{equation}
 is calm at $0$ for all $Y\!\in\!\Upsilon_0(0)$ (see Section \ref{sec2.1} for its definition), and developed a proximal DC algorithm for solving the  problem \eqref{EroneDNNcon-penalty} with a carefully selected $\rho$. 

 The global exact penalty result in \cite{Jiang21} is imperfect because it is unclear whether the calmness assumption on $\Upsilon_0$ at $0$ for all $Y\in\Upsilon_0(0)$ holds or not. In fact, to verify the calmness of a multifunction is not an easy task, since no convenient criterion is available. In addition, the solutions returned by the proximal DC method in \cite{Jiang21} indeed have high quality, but the required time is still unacceptable even for the medium-size $n\in\{31,\ldots,60\}$. The two aspects offer the motivation for this work. 
 \subsection{Main contributions}\label{sec1.2}

 Let $D\in\mathbb{S}^p$ with $D^{ii}=ee^{\top}\!-I$ and $D^{ij}=I$ for all $i\ne j\in[n]$, and let $\Lambda\in\mathbb{S}^p$ with $\Lambda^{ii}=0$ and $\Lambda^{ij}=I$ for all $i\ne j\in[n]$. Then, it is not difficult to check that for any $Y\!\in\mathbb{S}^p\cap\mathbb{R}_{+}^{p\times p}$, 
 \begin{subequations}
 \begin{align}\label{imply-set}
 \langle D,Y\rangle =0\ \Longrightarrow\ Y\in\mathbb{L}:=\big\{Y\in\mathbb{S}^{p}\ |\ \Lambda\circ Y=0\big\},\qquad\\
 \label{equiv-set}
 \langle I,Y^{ij}\rangle=\varpi_{ij}\ \ \forall i,j\in[n],\,{\textstyle\sum_{t=1}^n}Y^{tt}=I\Longleftrightarrow\qquad\qquad\\
 \langle D,Y\rangle\!=0,\,\langle I,Y^{ii}\rangle=1,\,{\textstyle\sum_{t=1}^n}(Y^{tt})_{ii}=1\ \ \forall i\in[n]. \nonumber
 \end{align} 
 \end{subequations} 
 By the expression of $\Delta_0$, the equivalence \eqref{equiv-set} implies that \eqref{EroneDNNcon} is equivalent to 
 \begin{align}\label{EroneDNNcon1}
 &\min_{Y\in\mathbb{S}^p}\ \langle G,Y\rangle\nonumber\\
 &\ {\rm s.t.}\ \ \langle I,Y^{ii}\rangle=1,\, {\textstyle\sum_{t=1}^n}{(Y^{tt})_{ii}}=1\ \ {\rm for}\  i\in[n],\nonumber\\
 &\qquad\ \langle D,Y\rangle=0,\ \langle ee^{\top},Y\rangle=p,\,Y\!\in K,\,{\rm tr}(Y)-\|Y\|=0,
 \end{align}
 whose feasible set, denoted by $\Gamma$, coincides with that of  \eqref{EroneDNNcon}. Define the sets
 \begin{subequations}
 \begin{align}\label{Tset-def}
 T:=\big\{Y\in\mathbb{S}^{p}\ |\ \langle I,Y^{ii}\rangle=1\ {\rm for}\ i\in[n]\big\},\qquad\qquad\qquad\qquad\\
 \Omega:=\big\{Y\in\mathbb{S}^p\ |\ \langle D,Y\rangle=0,\langle ee^{\top},Y\rangle=p,{\textstyle\sum_{t=1}^n}(Y^{tt})_{ii}=1\ {\rm for}\ i\in[n]\big\}.
 \label{Omega-def}
 \end{align}
 \end{subequations}
 Together with the definition of $\Delta_0$ and the above \eqref{imply-set}-\eqref{equiv-set}, it follows  
 \begin{equation}\label{imply1-set}
  \Delta_0=T\cap\Omega\cap K\subset T\cap\mathbb{L}\cap K,
 \end{equation}
 so $\Gamma=\Delta_0\cap\mathcal{R}=T\cap\Omega\cap K\cap\mathcal{R}\subset T\cap\mathbb{L}\cap K\cap\mathcal{R}$. In fact, from Lemma \ref{Gamma-set} later, 
 \begin{equation}\label{Eset-Gamma}
 \Gamma=  \Delta_0\cap\mathcal{R}=T\cap\Omega\cap K\cap\mathcal{R}=T\cap\mathbb{L}\cap K\cap\mathcal{R}.
 \end{equation}
 
 Our first contribution is to establish the locally and globally Lipschitzian error bounds for the set $\Gamma$ by its expressions $(T\cap\mathbb{L}\cap K)\cap\mathcal{R}$, $(T\cap\Omega\cap K)\cap\mathcal{R}$ and $\Delta_0\cap\mathcal{R}$, and employ the global error bounds to derive the locally upper Lipschitz property of 
 \begin{subequations}
 \begin{align}\label{Upsilon-map}
 \Upsilon(\tau)&:=\Big\{Y\in T\cap\mathbb{L}\cap K\ |\ {\rm tr}(Y)-\|Y\|=\tau\Big\}\quad{\rm for}\ \tau\in\mathbb{R},\\
 \label{Upsilon1-map}
  \Upsilon_{\!1}(\tau)&:=\Big\{Y\in T\cap\Omega\cap K\ |\ {\rm tr}(Y)-\|Y\|=\tau\Big\}\quad{\rm for}\ \tau\in\mathbb{R},
 \end{align}    
 \end{subequations} 
 and that of the previous $\Upsilon_0$ at $\tau=0$. Then, the calmness assumption required in \cite{Jiang21} for the mapping $\Upsilon_0$ holds automatically. It is worth emphasizing that many criteria were proposed for identifying calmness of a multifunction or metric subregularity of its inverse mapping (see, e.g., \cite{Henrion05,Zheng10,Gfrerer11,BaiYe19}), but they are not workable for checking the calmness of the multifunctions $\Upsilon,\Upsilon_1$ or $\Upsilon_0$ due to their complicated structure. 
 
 Let $f:\mathbb{S}^{p}\to\mathbb{R}$ be a function that is locally Lipschitz relative to the set $T\cap\mathbb{L}\cap K$, and let $\Sigma$ signify one of the sets $T\cap \mathbb{L}\cap K$, $T\cap \Omega\cap K$ and $\Delta_0$. Denote by $\interleave \cdot\interleave$ either the spectral norm $\|\cdot\|$ or the Frobenius norm $\|\cdot\|_F$. Consider the problem
 \begin{equation}\label{gDNN-mpec}
 \min_{Y\in\mathbb{S}^p}\Big\{f(Y)\ \ {\rm s.t.}\ \ {\rm tr}(Y)-\interleave Y\interleave=0,\,Y\in\Sigma\Big\}.
 \end{equation}
 The second contribution is to prove that, when a perturbation is imposed on the DC constraint ${\rm tr}(Y)-\interleave Y\interleave=0$, the problem \eqref{gDNN-mpec} is partially calm on the set of local optimal solutions. Along with the compactness of $\Sigma$, the following penalty problem 
 \begin{align}\label{gDNN-penalty}
 &\min_{Y\in\mathbb{S}^p}\Big\{f(Y)+\rho({\rm tr}(Y)-\interleave Y\interleave)\ \ {\rm s.t.}\ \ Y\in\Sigma\Big\}
 \end{align}
 is proved to be a global exact penalty of \eqref{gDNN-mpec}. This recovers the exact penalty result in \cite{Jiang21} for \eqref{EroneDNNcon-penalty} by removing the calmness assumption on the mapping $\Upsilon_0$.  

 When applying a convex relaxation approach to solve a single penalty problem \eqref{gDNN-penalty}, one eigenvalue decomposition for a matrix in $\mathbb{S}^p$ is at least required at each iteration, which forms the major computational bottleneck and restricts their scalability to large-scale problems. Motivated by the recent renewed interest in the BM factorization method \cite{Burer01,Burer03} for low-rank optimization problems (see, e.g., \cite{Lee2023,Tang2023,TaoQianPan22}), it is natural to develop relaxation methods by the BM factorization. Let $\Xi$ represent one of the compact sets $T\cap\mathbb{L}\cap\mathbb{R}_{+}^{p\times p},T\cap\Omega\cap\mathbb{R}_{+}^{p\times p}$ and $\Delta_0'$, where $\Delta_0'$ is the set obtained by replacing $K$ in $\Delta_0$ with $\mathbb{R}_{+}^{p\times p}$. Choose an appropriate $m\in[p]$. The BM factorization of the rank-one DNN constrained problem \eqref{gDNN-mpec} takes the form of 
 \begin{equation}\label{gfac-mpec}
 \min_{V\in\mathbb{R}^{m\times p}}\Big\{f(V^{\top}V)\ \ {\rm s.t.}\ \ \|V\|_F^2-\interleave V^{\top}V\interleave =0,\,V^{\top}V\in\Xi\Big\},
 \end{equation}
 where $V=[V_1\ V_2\ \cdots\ V_n]$ with $V_i\in\mathbb{R}^{m\times n}$ for each $i\in[n]$. The third contribution is to show that, when a perturbation is imposed on the DC constraint $\|V\|_F^2-\interleave V^{\top}V\interleave=0$, the problem \eqref{gfac-mpec} is partially calm on the set of its local optimal solutions. Together with the compactness of the set $\{V\in\mathbb{R}^{m\times p}\,|\,V^{\top}V\in\Xi\}$, the following problem
 \begin{equation}\label{gfac-penalty}
  \min_{V\in\mathbb{R}^{m\times p}}\Big\{f(V^{\top}V)
	+\rho(\|V\|_F^2-\interleave V^{\top}V\interleave)\ \ {\rm s.t.}\ \ V^{\top}V\in\Xi\Big\}
 \end{equation}
 is proved to be a global exact penalty for the nonconvex problem \eqref{gfac-mpec}. 
  
 Finally, with  \eqref{gfac-penalty} for $\Xi=T\cap\Omega\cap\mathbb{R}_{+}^{p\times p}$, we propose a relaxation algorithm to seek a rank-one approximate feasible point so as to illustrate the application of these exact penalties in developing such methods. This approach yields a satisfactory rank-one approximate feasible solution by searching for approximate stationary points of a finite number of penalty subproblems with an ALM (EPalm, for short). Its efficiency is validated by comparing the quality of the generated solutions with that of those given by Gurobi within the same running time for \textbf{132} benchmark instances, including $\textbf{122}$ QAPLIB instances \cite{QAPLIB} and $\textbf{10}$ drexxx instances \cite{Drezner15}. Numerical results show that EPalm is significantly superior to Gurobi for the QAPLIB instances with $n\le 60$ and the drexxx instances, and comparable with the latter for the QAPLIB instances with $n>60$, in terms of the relative gap between the generated objective value and the known best one and the number of examples with better objective values. For the $\textbf{10}$ drexxx instances, EPalm with a trivial rounding produces the known best values for all of them, while Gurobi yields the relative gaps more than $50\%$ for $\textbf{6}$ instances. For the $\textbf{108}$ QAPLIB instances with $n\le 60$, there are $\textbf{68}$ and $\textbf{9}$ ones by EPalm respectively better and worse than those by Gurobi. 
 For the $\textbf{14}$ QAPLIB instances with $n>60$, there are $\textbf{5}$ ones by EPalm respectively with better and worse objective values, and the maximum relative gap $8.24\%$ is greater than that of Gurobi $5.93\%$. When comparing the objective values by EPalm with the upper bounds by rPRSM of \cite{Graham22}, for the total $\textbf{84}$ test instances reported in \cite[Tables 1-3]{Graham22}, EPalm returns the tighter upper bounds for $\textbf{53}$ instances, and the worse upper bounds for $\textbf{3}$ ones. 
\subsection{Notation}\label{sec1.3}

 Throughout this paper, $\mathbb{S}^{p}$ represents the set of all $p\times p$ real symmetric matrices, $\mathbb{S}_{+}^{p}$ denotes the cone consisting of all $p\times p$ positive semidefinite matrices, $\mathbb{R}_{+}^{m\times p}$ signifies the cone of all $m\times p$ nonnegative matrices, and $\mathbb{O}^{p}$ denotes the set of all $p\times p$ orthonormal matrices. The notation $I$ and $e$ respectively denotes an identity matrix and a vector of all ones, whose dimensions are known from the context. For an integer $k\ge 1$, write $[k]:=\{1,\ldots,k\}$ and $[k]_{+}:=\{0\}\cup[k]$. For any $Z\in\mathbb{R}^{n\times q}$, $\|Z\|_*$ and $\|Z\|$ denote the nuclear norm and spectral norm of $Z$, respectively, $Z_j$ represents the $j$th column of $Z$, and $\mathbb{B}(Z,\varepsilon)$ denotes the closed ball on the Frobenius norm centered at $Z$ with radius $\varepsilon$. The notation ${\rm vec}(Z)$ for a matrix signifies a column vector whose entries come from $Z$ by stacking up columns from the first to the last column on top of each other, and the operator ${\rm mat}$ denotes the inverse of ${\rm vec}$, i.e., ${\rm mat}({\rm vec}(Z))=Z$. For any $Y\in\mathbb{S}^p$, let $\mathbb{O}^p(Y):=\big\{P\in\mathbb{O}^{p}\ |\ Y=P{\rm Diag}(\lambda(Y))P^{\top}\big\}$ where $\lambda(Y)$ is the eigenvalue vector of $Y$ arranged in a nonincreasing order, and $Y^{ij}$ for $i,j\in[n]$ denotes the $(i,j)$-th block submatrix of $Y$ with dimension $n\times n$. For a closed set $\Delta\subset\mathbb{R}^{n\times q}$, $\chi_\Delta$ denotes the indicator function of $\Delta$, $\Pi_{\Delta}(\cdot)$ denotes the projection mapping on $\Delta$, and ${\rm dist}(X,\Delta)$ means the distance on the Frobenius norm from $X$ to $\Delta$. 
 
\section{Preliminaries}\label{sec2}

In this section, a hollow capital letter, say $\mathbb{X}$, signifies a finite-dimensional real Euclidean space equipped with the inner product $\langle\cdot,\cdot\rangle$ and its induced norm $\|\cdot\|$, and $\mathbb{B}_{\mathbb{X}}$ denotes the unit closed ball of $\mathbb{X}$ centered at the origin on $\|\cdot\|$. 
\subsection{Calmness and subregularity}\label{sec2.1}

 The notion of calmness of a multifunction was first introduced in \cite{YeYe97} under the term ``pseudo upper-Lipschitz continuity'' owing to the fact that it is a combination of Aubin's pseudo-Lipschitz continuity and Robinson's local upper-Lipschitz continuity \cite{Robinson81}, and the term ``calmness'' was later coined in \cite{RW98}. A multifunction $\mathcal{S}\!:\mathbb{Y}\rightrightarrows\mathbb{Z}$ is said to be calm at $\overline{y}$ for $\overline{z}\in\mathcal{S}(\overline{y})$ if there exists $\nu\ge0$ along with $\varepsilon>0$ and $\delta>0$ such that for all $y\in\mathbb{B}(\overline{y},\varepsilon)$, 
 \begin{equation}\label{calm-def}
 \mathcal{S}(y)\cap\mathbb{B}(\overline{z},\delta)
 \subset\mathcal{S}(\overline{y})+\nu\|y-\overline{y}\|\mathbb{B}_{\mathbb{Z}}.
 \end{equation}
 By \cite[Exercise 3H.4]{DR09}, the neighborhood restriction on $y$ in \eqref{calm-def} can be removed. As observed by Henrion and Outrata \cite{Henrion05}, the calmness of $\mathcal{S}$ at $\overline{y}$ for $\overline{z}\in\mathcal{S}(\overline{y})$ is equivalent to the (metric) subregularity of its inverse at $\overline{z}$ for $\overline{y}\in\mathcal{S}^{-1}(\overline{z})$. Subregularity was introduced by Ioffe in \cite{Ioffe79} (under a different name) as a constraint qualification related to equality constraints in nonsmooth optimization problems, and was later extended to generalized equations. Recall that a multifunction $\mathcal{F}\!:\mathbb{Z}\rightrightarrows\mathbb{Y}$ is called (metrically) subregular at $\overline{z}$ for $\overline{y}\in\mathcal{F}(\overline{z})$ if there exists $\kappa'\ge 0$ along with $\varepsilon>0$ such that
 \begin{equation*}
  {\rm dist}(z,\mathcal{F}^{-1}(\overline{y}))\le \kappa'{\rm dist}(\overline{y},\mathcal{F}(z))
 \quad{\rm for\ all}\ z\in\mathbb{B}(\overline{z},\varepsilon).
 \end{equation*}
 The calmness and subregularity have been studied by many authors under various names (see, e.g., \cite{Henrion05,Ioffe08,Gfrerer11,BaiYe19,Zheng10}
 and the references therein). From \cite{Robinson81}, $\mathcal{S}$ is said to be locally upper Lipschitzian at $\overline{y}$ with modulus $\gamma$ if there exists $\varepsilon>0$ such that 
 \begin{equation}\label{UL-def}
 \mathcal{S}(y)\subset\mathcal{S}(\overline{y})+\gamma\|y-\overline{y}\|\mathbb{B}_{\mathbb{Y}}\quad{\rm for\ all}\ y\in\mathbb{B}(\overline{y},\varepsilon).
 \end{equation}
 Clearly, the local upper Lipschitz of $\mathcal{S}$ at $\overline{y}$ implies its calmness at $\overline{y}$ for any $\overline{z}\in\mathcal{S}(\overline{y})$. 
\subsection{Partial calmness of optimization problems}\label{sec2.2}

 Let $\vartheta\!:\mathbb{Z}\to (-\infty,\infty]$ be a proper lower semicontinuous (lsc) function and  $h\!:\mathbb{Z}\to\mathbb{R}$ be a continuous function. This section focuses on the calmness of the following abstract optimization problem
 \begin{equation*}
 ({\rm MP})\qquad\min_{z\in\mathbb{Z}}\big\{\vartheta(z)\ \ {\rm s.t.}\ \ h(z)=0,\,z\in\Delta\big\}
 \end{equation*}
 when only $h(z)=0$ is perturbed, where $\Delta\subset\mathbb{Z}$ is a nonempty closed set. This plays a key role in achieving the local and global exact penalty induced by the constraint $h(z)=0$ (see \cite{YeZhu97,LiuBiPan18}). The calmness of an optimization problem at a solution point was first introduced by Clarke \cite{Clarke83}. Later, Ye and Zhu \cite{YeYe97,YeZhu97} extended it to the partial calmness at a solution point, whose formal definition is stated as follows.
 \begin{definition}\label{def-pcalm0}(see \cite[Definition 2.1]{YeZhu97}) Denote by $\mathcal{Z}^*$ the set of local optimal solutions of $({\rm MP})$, and define $\mathcal{S}(\tau):=\big\{z\in\Delta\ |\ h(z)=\tau\big\}$ for $\tau\in\mathbb{R}$. The $({\rm MP})$ is said to be partially calm at a point $z^*\in\mathcal{Z}^*$ if there exist $\varepsilon>0$ and $\mu>0$ such that 
 \[
   \vartheta(z)-\vartheta(z^*)+\mu |h(z)|\geq 0
 \]
 for all $\tau\in\mathbb{R}$ and all $z\in(z^*\!+\varepsilon\mathbb{B}_{\mathbb{Z}})\cap\mathcal{S}(\tau)$. If the problem $({\rm MP})$ is partially calm at every point of $\mathcal{Z}^*$, it is said to be partially calm on the set $\mathcal{Z}^*$.
 \end{definition}

 Definition \ref{def-pcalm0} has a little difference  from \cite[Definition 2.1]{YeZhu97}. Since the function $h$ is assumed to be continuous, by virtue of \cite[Remark 2.3]{YeZhu97}, the restriction on the size of perturbation $\tau$ is removed from Definition \ref{def-pcalm0}. The following lemma states that the partial calmness of $({\rm MP})$ on $\mathcal{Z}^*$ is implied by the calmness of $\mathcal{S}$ at $0$ for all $z\in\mathcal{Z}^*$. The proof is similar to that of \cite[Lemma 3.1]{YeYe97}, and we include it for completeness.
 \begin{lemma}\label{lemma-calm}
 If $\vartheta$ is locally Lipschitz continuous relative to $\Delta$ and the mapping $\mathcal{S}$ in Definition \ref{def-pcalm0} is calm at $0$ for any $z\in\mathcal{Z}^*$, then $({\rm MP})$ is partially calm on $\mathcal{Z}^*$.
 \end{lemma}
 \begin{proof}
 Pick any $z^*\in\mathcal{Z}^*$. Since $z^*\in\Delta$ and $\vartheta$ is locally Lipschitz relative to $\Delta$, there exist $\varepsilon'>0$ and $L_{\vartheta}>0$ such that for all $z',z''\in\mathbb{B}(z^*,\varepsilon')\cap\Delta$,
 \begin{equation}\label{temp-equa21}
 |\vartheta(z')-\vartheta(z'')|\le L_{\vartheta}\|z'-z''\|.
 \end{equation}
 Notice that $z^*$ is a local optimal solution of $({\rm MP})$. If necessary by shrinking $\varepsilon'$, 
 \begin{equation}\label{temp-equa22}
 \vartheta(z)\ge\vartheta(z^*)\quad{\rm for\ all}\ z\in\mathbb{B}(z^*,\varepsilon')\cap\mathcal{S}(0).
 \end{equation}
 Since the mapping $\mathcal{S}$ is calm at $0$ for $z^*$, there exist $\nu>0$ and $\delta'>0$ such that 
 \begin{equation*}
 \mathcal{S}(\omega)\cap\mathbb{B}(z^*,\delta')\subset\mathcal{S}(0)+\nu|\omega|\mathbb{B}_{\mathbb{Z}}\quad{\rm for\ all}\ \omega\in\mathbb{R}.
 \end{equation*}
 Set $\varepsilon=\frac{\min\{\varepsilon',\delta'\}}{2}$. Pick any $\tau\in\mathbb{R}$ and $z\in(z^*\!+\!\varepsilon\mathbb{B}_{\mathbb{Z}})\cap\mathcal{S}(\tau)$. Since $z\in\mathbb{B}(z^*,\delta')\cap\mathcal{S}(\tau)$, using the above inclusion with $\omega=\tau$ leads to ${\rm dist}(z,\mathcal{S}(0))\le\nu|\tau|=\nu|h(z)|$. From the closedness of $\mathcal{S}(0)$, there exists  $\widehat{z}\in\mathcal{S}(0)$ such that $\|z-\widehat{z}\|={\rm dist}(z,\mathcal{S}(0))\le\nu|h(z)|$. Note that $\|\widehat{z}-z^*\|\le 2\|z-z^*\|\le\varepsilon'$. Together with $\widehat{z}\in\mathcal{S}(0)$ and $z\in\Delta$, 
 \[
   \vartheta(z^*)\stackrel{\eqref{temp-equa22}}{\le} \vartheta(\widehat{z})=\vartheta(z)-\vartheta(z)+\vartheta(\widehat{z})
   \stackrel{\eqref{temp-equa21}}{\le} \vartheta(z)+L_{\vartheta}\|z-\widehat{z}\|\le \vartheta(z)+L_{\vartheta}\nu|h(z)|,
 \]
 so $({\rm MP})$ is partially calm at $z^*$. The conclusion follows the arbitrariness of $z^*\!\in\mathcal{Z}^*$.
 \end{proof}
\section{Lipschitzian error bounds}\label{sec3}

This section focuses on the following locally Lipschitzian error bound: to seek a constant $\kappa>0$ (depending only on $n$) such that for each $\overline{Y}\!\in\Gamma$, there exists $\varepsilon>0$ such that for all $Y\!\in\mathbb{B}(\overline{Y},\varepsilon)$,
\begin{equation}\label{aim-error}
{\rm dist}(Y,\Gamma)
 \le \kappa\big[{\rm dist}(Y,T)+{\rm dist}(Y,\mathbb{L})+4{\rm dist}(Y,K)+{\rm dist}(Y,\mathcal{R})\big],
\end{equation}
which implies the locally and globally Lipschitzian error bounds of $\Gamma$ by its expressions $(T\cap \mathbb{L}\cap K)\cap\mathcal{R},(T\cap \Omega\cap K)\cap\mathcal{R}$ and $\Delta_0\cap\mathcal{R}$. First, we give a technical lemma on $\Gamma$. 
 \begin{lemma}\label{Gamma-set}
 For the feasible set $\Gamma$ of the problem \eqref{EroneDNNcon1},  it holds that
 \[	
  \Gamma=\widehat{\Gamma}:=\big\{{\rm vec}(X){\rm vec}(X)^{\top}\, |\, X^{\top}X=I,X\in\mathbb{R}_{+}^{n\times n}\big\}=T\cap\mathbb{L}\cap K\cap\mathcal{R}.
 \]
 \end{lemma}
 \begin{proof}
 Pick any $Y\in\Gamma$. Recall that $\Gamma=T\cap\Omega\cap K\cap\mathcal{R}$. From $Y\in K\cap\mathcal{R}$, there exists $x\in\mathbb{R}_{+}^p$ such that $Y=xx^{\top}$, which along with $Y\in T\cap\Omega$ means that $X={\rm mat}(x)\in\mathbb{R}^{n\times n}$ satisfies $X^{\top}X=I$. Then, $Y\in\widehat{\Gamma}$ and $\Gamma\subset\widehat{\Gamma}$ holds. For the converse inclusion, pick any $Y\in\widehat{\Gamma}$. Then, $Y={\rm vec}(X){\rm vec}(X)^{\top}$ for some $X\in\mathbb{R}_{+}^{n\times n}$ with $X^{\top}X=I$. Clearly, $Y\in K\cap\mathcal{R}$. Notice that each row and each column of $X$ has only one nonzero entry $1$. It is easy to check $Y\in T\cap\Omega$, so $Y\in\Gamma$ and $\Gamma\supset\widehat{\Gamma}$ follows. The first equality holds. The inclusion \eqref{imply1-set} implies $\Gamma=T\cap\Omega\cap K\cap\mathcal{R}\subset T\cap\mathbb{L}\cap K\cap\mathcal{R}$. Then, for the second equality, it suffices to prove that $T\cap\mathbb{L}\cap K\cap\mathcal{R}\subset\widehat{\Gamma}$. Pick any $Y\in T\cap\mathbb{L}\cap K\cap\mathcal{R}$. From $Y\in K\cap\mathcal{R}$, we infer that $Y=xx^{\top}$ with $x\in\mathbb{R}_{+}^p$. Let $X={\rm mat}(x)$. Then, $X\in\mathbb{R}_{+}^{n\times n}$ and $Y={\rm vec}(X){\rm vec}(X)^{\top}$. Together with $Y\in T\cap\mathbb{L}$, we obtain $X^{\top}X=I$, so $Y\in\widehat{\Gamma}$. The desired inclusion holds. 
\end{proof}
 
 To achieve the local error bound in \eqref{aim-error}, we next establish three propositions. The first one states that for any $\overline{Y}\in\Gamma$ and any $Y$ close to $\overline{Y}$, the distance of $Y$ from the set $\Gamma$ is locally upper bounded by the sum of ${\rm dist}(Y, T)$ and ${\rm dist}(Y,\mathbb{L}\cap K\cap\mathcal{R})$. Among others, the discreteness of $\Gamma$ and the structure of $T$ play a crucial role. 
 \begin{proposition}\label{dist-prop1}
 Fix any $\overline{Y}\in\Gamma$. There is $\delta>0$ such that for all $Y\in\mathbb{B}(\overline{Y},\delta)$,
 \begin{equation*}
 {\rm dist}(Y,\Gamma)\!\le\!\big(1\!+\!\sqrt{n(4n\!+\!1)}\big)\big[{\rm dist}(Y,T)+{\rm dist}(Y,\mathbb{L}\cap K \cap\mathcal{R})\big].
 \end{equation*}
\end{proposition}
\begin{proof}
 From $\overline{Y}\in\Gamma$ and the first equality of Lemma \ref{Gamma-set}, there exists $\overline{Z}\in\mathbb{R}_{+}^{n\times n}$ with $\overline{Z}^{\top}\overline{Z}=I$ such that $\overline{Y}={\rm vec}(\overline{Z}){\rm vec}(\overline{Z})^{\top}$. Since each row and each column of $\overline{Z}$ has only a nonzero entry $1$, every $\overline{Y}^{ii}\in\mathbb{S}^n$ for $i\in[n]$ has only one nonzero entry $1$ in its diagonal, and every $\overline{Y}^{ij}\in\mathbb{R}^{n\times n}$ for $i\ne j\in[n]$ has only one nonzero entry $1$. The first equality of Lemma \ref{Gamma-set} implies that the entries of every $Y\in\Gamma$ all belong to $\{0,1\}$. By the discreteness of $\Gamma$ and $\overline{Y}\!\in\Gamma$, there exists $\delta'\in(0,\frac{1}{2}]$ such that
 \begin{equation}\label{equa-MF}
 {\rm dist}(Y,\Gamma)=\|Y-\overline{Y}\|_F\quad\ \forall\, Y\in\mathbb{B}(\overline{Y},\delta').
 \end{equation} 
 
 Next we prove ${\rm dist}(Y,\Gamma)\le\!\sqrt{n(4n\!+\!1)}\,{\rm dist}(Y,T)$ for all $Y\in\mathbb{B}(\overline{Y},\delta')\cap[\mathbb{L}\cap K\cap\mathcal{R}]$. 
 Fix any $Y\in\mathbb{B}(\overline{Y},\delta')\cap[\mathbb{L}\cap K\cap\mathcal{R}]$. From $Y\in K\cap\mathcal{R}$, there necessarily exists a matrix $Z\in\mathbb{R}_{+}^{n\times n}$ such that $Y={\rm vec}(Z){\rm vec}(Z)^{\top}$ (if not, there is a matrix $Z'\in\mathbb{R}^{n\times n}$ such that $Y={\rm vec}(Z'){\rm vec}(Z')^{\top}$ with $[{\rm vec}(Z')]_j<0$ for some $j\in[p]$ and there exists an index $[p]\ni i\ne j$ such that $[{\rm vec}(Z')]_i>0$ (if such an index $i$ does not exist, $-Z'\in\mathbb{R}_{+}^{n\times n}$ is such that  
 $Y={\rm vec}(-Z'){\rm vec}(-Z')^{\top}$ and satisfies the desired one). Then, we have $Y_{ij}<0$, a contradiction to the fact that $Y\in K$). Along with $Y\in\mathbb{L}$ and the definition of $E$, we get $\langle Z_{i},Z_{j}\rangle=0$ for all $i\neq j\in[n]$. Hence, each row of $Z$ has at most a nonzero entry and $Z$ has at most $n$ nonzero entries. Recall that every $\overline{Y}^{ii}$ has only one nonzero entry $1$ appearing in its diagonal. From $Y\in\mathbb{B}(\overline{Y},\delta')$, we infer that every $Y^{ii}$ has at least one nonzero entry in its diagonal, which along with ${\rm diag}(Y^{ii})=Z_{i}\circ Z_{i}$ for each $i\in[n]$ means that $Z_{i}$ has at least one nonzero entry. Since $Z$ has at most $n$ nonzero entries, each column of $Z$ has only a nonzero entry. Consequently, for each $i\in[n]$, ${\rm diag}(Y^{ii})$ has only a nonzero entry. For each $i\in[n]$, let $\overline{t}_i,t_i\in[n]$ be such that  $[\overline{Y}^{ii}]_{\overline{t}_i\overline{t}_i}=1$ and $[Y^{ii}]_{t_it_i}>0$. Then, $t_i=\overline{t}_i$ and
 \begin{equation}\label{equa-Tset}
  [{\rm dist}(Y,T)]^2=\frac{1}{n^2}\sum_{i=1}^n\big\|(1-\langle I,Y^{ii}\rangle)I\big\|_F^2=\frac{1}{n}\sum_{i=1}^n|1-(Y^{ii})_{t_it_i}|^2,
 \end{equation}
 where the first equality is due to the definition of $T$ in \eqref{Tset-def}. While from \eqref{equa-MF}, $Y={\rm vec}(Z){\rm vec}(Z)^{\top}$ and $\overline{Y}={\rm vec}(\overline{Z}){\rm vec}(\overline{Z})^{\top}$, it follows 
 \begin{align*}
  [{\rm dist}(Y,\Gamma)]^2
 &=\big\|{\rm vec}(Z){\rm vec}(Z)^{\top}-{\rm vec}(\overline{Z}){\rm vec}(\overline{Z})^{\top}\big\|_F^2\\
 &\le(\|Z\|_F+\|\overline{Z}\|_F)^2\|Z-\overline{Z}\|_F^2\\
 &=(\|Z\|_F+\|\overline{Z}\|_F)^2\,\textstyle{\sum_{i=1}^n}\big|\sqrt{(Y^{ii})_{t_it_i}}-1\big|^2.
 \end{align*}
 Note that $\|Z\|_F^2=\|{\rm vec}(Z)\|^2=\sum_{i=1}^n\|{\rm diag}(Y^{ii})\|_1=\|{\rm diag}(Y)\|_1=\langle I,Y\rangle=\|Y\|_F$, where the last two equalities is due to $Y\in K\cap\mathcal{R}$. Along with $\|Y\|_F\le\|\overline{Y}\|_F+\delta'\le n+{1}/{2}$, we have $\|Z\|_F^2\le n+{1}/{2}$, so  $(\|Z\|_F+\|\overline{Z}\|_F)^2\le 2(\|Z\|_F^2+\|\overline{Z}\|_F^2)\le 2(2n+1/2)$. Then, from the above inequality and $|1-\!\sqrt{(Y^{ii})_{t_it_i}}|=\frac{|1-(Y^{ii})_{t_it_i}|}{1+\sqrt{(Y^{ii})_{t_it_i}}}\le |1-(Y^{ii})_{t_it_i}|$, 
 \begin{displaymath}
  [{\rm dist}(Y,\Gamma)]^2
  \le (4n\!+1)\sum_{i=1}^n|1-(Y^{ii})_{t_it_i}|^2.
 \end{displaymath}
 Combining this inequality with \eqref{equa-Tset} shows that the desired inequality holds. 
 
 Now fix any $Y\in\mathbb{B}(\overline{Y},\delta)$ with $\delta=\delta'/2$. Pick any $\widetilde{Y}\in\Pi_{\mathbb{L}\cap K\cap\mathcal{R}}(Y)$. Obviously, $\|\widetilde{Y}-\overline{Y}\|_F\le 2\|\overline{Y}-Y\|_F\le\delta'$. Then, from the above arguments, it follows ${\rm dist}(\widetilde{Y},\Gamma)\le \sqrt{n(4n\!+\!1)}\,{\rm dist}(\widetilde{Y},T)$. By leveraging this inequality, we immediately obtain
 \begin{align*}
  {\rm dist}(Y,\Gamma)
  &\le{\rm dist}(\widetilde{Y},\Gamma)+\|Y\!-\!\widetilde{Y}\|_F\le\sqrt{n(4n\!+\!1)}\,{\rm dist}(\widetilde{Y},T)+\|Y\!-\!\widetilde{Y}\|_F\nonumber\\
  &\le\sqrt{n(4n\!+\!1)}\,{\rm dist}(Y,T)+\big(1\!+\!\sqrt{n(4n\!+\!1)}\big)\|Y\!-\!\widetilde{Y}\|_F\\
  &\le\big(1\!+\!\sqrt{n(4n\!+\!1)}\big)\,\big[{\rm dist}(Y,T)+{\rm dist}(Y,\mathbb{L}\cap K\cap\mathcal{R})\big].
 \end{align*}
 This shows that the desired conclusion holds. Thus, we complete the proof.
\end{proof}

For the term ${\rm dist}(Y,\mathbb{L}\cap K \cap\mathcal{R})$ in Proposition \ref{dist-prop1}, the following proposition shows that it can be locally controlled by the sum of ${\rm dist}(Y,\mathbb{L})$ and ${\rm dist}(Y,K\cap\mathcal{R})$.
 \begin{proposition}\label{dist-prop2}
 Fix any $\overline{Y}\in\Gamma$. There is $\delta>0$ such that for all $Y\in\mathbb{B}(\overline{Y},\delta)$,
 \begin{equation*}
 {\rm dist}(Y,\mathbb{L}\cap\mathcal{R}\cap K)
 \le(1\!+\!3\sqrt{2}n)\big[{\rm dist}(Y,\mathbb{L})+{\rm dist}(Y,K\cap\mathcal{R})\big].
 \end{equation*}
\end{proposition}
\begin{proof}
 From $\overline{Y}\in\Gamma$ and Lemma \ref{Gamma-set},  $\Lambda\circ\overline{Y}=0$ and there exists $\overline{x}\in\{0,1\}^{p}$ such that $\overline{Y}=\overline{x}\overline{x}^{\top}$. For each $i\in[n]$, since $\sum_{t=1}^n(\overline{Y}^{tt})_{ii}=1$, it holds $\sum_{k=0}^{n-1}\overline{x}_{i+kn}=1$. We claim that there exists $\delta'>0$ such that for any $yy^{\top}\in\mathbb{B}(\overline{Y},\delta')$ with $y\in\mathbb{R}_+^p$,
 \begin{equation}\label{ineq-yxbar}
  \|y-\overline{x}\|_{\infty}\le 1/2.
 \end{equation}
 If not, there exists a sequence $\{y^k(y^k)^{\top}\}_{k\in\mathbb{N}}$ with $y^k\in\mathbb{R}_+^p$ and $y^k(y^k)^{\top}\in\mathbb{B}(\overline{Y},\frac{1}{k})$ such that $|y^k_{i_k}-\overline{x}_{i_k}|>\frac{1}{2}$ for some  $i_k\in\mathbb{N}$. Since $y_{i_k}^k\ge0$ and $\overline{x}_{i_k}\ge0$, we have $y_{i_k}^k+\overline{x}_{i_k}>\frac{1}{2}$. Then, $\|y^k(y^k)^{\top}-\overline{x}\overline{x}^{\top}\|_F 	\ge|(y_{i_k}^k)^2-\overline{x}_{i_k}^2|=(y_{i_k}^k+\overline{x}_{i_k})|y_{i_k}^k-\overline{x}_{i_k}|>\frac{1}{4}$, a contradiction to $y^k(y^k)^{\top}\in\mathbb{B}(\overline{Y},\frac{1}{k})$ for all $k$. The claimed inequality \eqref{ineq-yxbar} holds.
 
 Let $\delta={\delta'}/{2}$. Fix any $Y\in\mathbb{B}(\overline{Y},\delta)$. Let $u\in\mathbb{R}_+^p$ be such that $uu^{\top}\!\in\Pi_{K\cap\mathcal{R}}(Y)$. 

 \noindent
 {\bf Step 1: to construct a vector $\overline{u}\in\mathbb{R}_{+}^p$ such that $\overline{u}\overline{u}^{\top}\in \mathbb{L}\cap K\cap\mathcal{R}$, so that} 
 \begin{align*}
 {\rm dist}(Y,\mathbb{L}\cap K\cap\mathcal{R})
 &\le\|Y-uu^{\top}\|_F+{\rm dist}(uu^{\top},\mathbb{L}\cap K\cap\mathcal{R})\nonumber\\
 &\le\|Y-uu^{\top}\|_F+\|uu^{\top}\!-\overline{u}\overline{u}^{\top}\|_F.
 \end{align*}
 To this end, write $\mathcal{I}(u,\overline{x}):=\big\{i\in[p]\ |\ u_i\neq0,\,\overline{x}_i=0\big\}$. Let $\overline{u}\in\mathbb{R}_+^p$ with $\overline{u}_i=0$ for $i\in\mathcal{I}(u,\overline{x})$ and $\overline{u}_i=u_i$ for $i\notin \mathcal{I}(u,\overline{x})$. Obviously, ${\rm supp}(\overline{u})\subset{\rm supp}(\overline{x})$. Along with $\Lambda\circ(\overline{x}\overline{x}^{\top})=\Lambda\circ\overline{Y}=0$, we have $\Lambda\circ(\overline{u}\overline{u}^{\top})=0$, so $\overline{u}\overline{u}^{\top}\in \mathbb{L}\cap K\cap\mathcal{R}$. 
 
 \noindent
 {\bf Step 2: to bound $\|uu^{\top}\!-\overline{u}\overline{u}^{\top}\|_F^2$.} For each $i\notin\mathcal{I}(u,\overline{x})$, since $(uu^{\top}\!-\overline{u}\overline{u}^{\top})_{ij}=u_iu_j$ for $j\in\mathcal{I}(u,\overline{x})$ and $(uu^{\top}\!-\overline{u}\overline{u}^{\top})_{ij}=0$ for $j\notin\mathcal{I}(u,\overline{x})$, it is immediate to obtain   
 \begin{equation}\label{temp-uubar}
  \|uu^{\top}\!-\overline{u}\overline{u}^{\top}\|_F^2\le 2\sum_{i\in\mathcal{I}(u,\overline{x})}\sum_{j=1}^p u_i^2u_j^2.
 \end{equation}
 Next we claim that for each $i\in\mathcal{I}(u,\overline{x})$, there exists an index $j_i\in[p]$ such that $\Lambda_{ij_i}=1$ and $\overline{x}_{j_i}=1$, which will be used to bound the sum on the right hand side of \eqref{temp-uubar}. If not, since each entry of $\Lambda$ and $\overline{x}$ belong to $\{0,1\}$, there exists an index $i_0\in\mathcal{I}(u,\overline{x})$ such that $\Lambda_{i_0j}\overline{x}_j=0$ for all $j\in[p]$. For the index $i_0$, there exist $m_0\in[n]$ and $t_0\in[n\!-\!1]_{+}$ such that $i_0=m_0+t_0n$ and $\overline{x}_{m_0+t_0n}=0$. Recalling that $\sum_{k=0}^{n-1}\overline{x}_{i+kn}=1$ for each $i\in[n]$, we have $\sum_{k=0}^{n-1}\overline{x}_{m_0+kn}=1$. While the definition of $\Lambda$ implies $\Lambda_{m_0+t_0n,m_0+kn}=1$ for each $k\in[n\!-\!1]_{+}\backslash\{t_0\}$. Then, it follows
 \begin{align*}
 \sum_{k\neq t_0,k=0}^{n-1}\Lambda_{m_0+t_0n,m_0+kn}\overline{x}_{m_0+kn}
 =\sum_{k\neq t_0,k=0}^{n-1}\overline{x}_{m_0+kn}=\sum_{k=0}^{n-1}\overline{x}_{m_0+kn}=1.
 \end{align*} 
 This, along with $m_0+t_0n=i_0$, implies that $\Lambda_{i_0\overline{j}}=1$ for some $\overline{j}\in[p]$, a contradiction to the fact that $\Lambda_{i_0j}\overline{x}_j=0$ for all $j\in[p]$. Thus, the claimed conclusion holds. Note that $\|uu^{\top}\!-\overline{Y}\|_F\le\|uu^{\top}\!-Y\|_F+\|Y-\overline{Y}\|_F\le 2\delta=\delta'$, so the above \eqref{ineq-yxbar} implies that $\|u-\overline{x}\|_{\infty}\le 1/2$. Recalling that $\overline{x}\in\{0,1\}^{p}$, we have $u_j\le 3/2$ for all $j\in[p]$. In addition, for each $i\in\mathcal{I}(u,\overline{x})$, 
 $u_{j_i}\ge\overline{x}_{j_i}-1/2=1/2$. 
 Thus, it holds 
 \[
  \sum_{i\in \mathcal{I}(u,\overline{x})}\sum_{j=1}^pu_i^2u_j^2\le ({9p}/{4})\!\sum_{i\in\mathcal{I}(u,\overline{x})}u_i^2\le 9p\sum_{i\in\mathcal{I}(u,\overline{x})}u_i^2u_{j_i}^2,
 \]
 Together with the above \eqref{temp-uubar}, the term $\|uu^{\top}\!-\overline{u}\overline{u}^{\top}\|_F^2$ is upper bounded as below
 \begin{equation}\label{temp2-uubar}
   \|uu^{\top}\!-\overline{u}\overline{u}^{\top}\|_F^2\le 18p\sum_{i\in\mathcal{I}(u,\overline{x})}u_i^2u_{j_i}^2.
 \end{equation}
 
 Let $J:=\big\{(i,j)\in[p]\times [p]\,|\,\Lambda_{ij}=1\big\}$. From Step 2, $\Lambda_{ij_i}=1$ for each $i\in\mathcal{I}(u,\overline{x})$, so $\sum_{i\in\mathcal{I}(u,\overline{x})}u_i^2u_{j_i}^2\le\sum_{(i,j)\in J}u_i^2u_{j}^2$. Noting that the projection of any $Y\in\mathbb{S}^p$ onto $\mathbb{L}$ satisfies $[\Pi_{\mathbb{L}}(Y)]_{ij}=0$ for $(i,j)\in J$, we have $[{\rm dist}(Y,\mathbb{L})]^2=\sum_{(i,j)\in J}Y_{ij}^2$. Then  
 \[
  [{\rm dist}(uu^{\top},\mathbb{L})]^2=\sum_{(i,j)\in J}u_i^2u_j^2\ge\sum_{i\in\mathcal{I}(u,\overline{x})}u_i^2u_{j_i}^2. 
 \]
 Now combining this inequality with \eqref{temp2-uubar} and the one in Step 1, it follows   
 \begin{align*}
  {\rm dist}(Y,\mathbb{L}\cap K\cap\mathcal{R})
  &\le\|Y\!-uu^{\top}\|_F+\sqrt{18p}\sqrt{{\textstyle\sum_{i\in\mathcal{I}(u,\overline{x})}}u_i^2u_{j_i}^2}\nonumber\\
  &\le\|Y\!-uu^{\top}\|_F+\sqrt{18p}\,{\rm dist}(uu^{\top},\mathbb{L})\nonumber\\
  &\le(1+\!\sqrt{18p}){\rm dist}(Y,K\cap\mathcal{R})+\sqrt{18p}\,{\rm dist}(Y,\mathbb{L}),
 \end{align*}
 where the third one is due to $uu^{\top}\!\in\!\Pi_{K\cap\mathcal{R}}(Y)$. The conclusion follows $p=n^2$ and the arbitrariness of $Y\in\mathbb{B}(\overline{Y},\delta)$. The proof is then completed.
 \end{proof}

 For the term ${\rm dist}(Y,K\cap\mathcal{R})$ in Proposition \ref{dist-prop2}, we have the following conclusion. 
\begin{proposition}\label{dist-KpR}
 For any $Y\in\mathbb{S}^p\cap\mathbb{R}_{+}^{p\times p}$, the following relations hold 
 \[
 {\rm dist}(Y,K\cap\mathcal{R})={\rm dist}(Y,\mathbb{S}_{+}^p\cap\mathcal{R})={\rm dist}(Y,\mathcal{R}).
 \]
\end{proposition}
\begin{proof}
 Fix any $Y\in\mathbb{S}^p\cap\mathbb{R}_{+}^{p\times p}$. According to \cite[Theorem 8.3.1]{Horn87}, the matrix $Y$ has an eigenvalue decomposition $U^{\top}{\rm Diag}(\lambda(Y))U$ with $U_1\in\mathbb{R}_{+}^p\backslash\{0\}$ corresponding to $\|Y\|$. Notice that ${\rm dist}(Y,\mathbb{S}_{+}^p\cap\mathcal{R})\ge{\rm dist}(Y,\mathcal{R})$, and $\|Y\|U_1U_1^{\top}\in\mathbb{S}_{+}^p\cap\mathcal{R}$ is such that $\|Y-\|Y\|U_1U_1^{\top}\|_F={\rm dist}(Y,\mathcal{R})$. Therefore, ${\rm dist}(Y,\mathbb{S}_{+}^p\cap\mathcal{R})={\rm dist}(Y,\mathcal{R})$. In addition, since ${\rm dist}(Y,K\cap\mathcal{R})\ge{\rm dist}(Y,\mathbb{S}_{+}^p\cap\mathcal{R})$ and $\|Y\|U_1U_1^{\top}\in K\cap\mathcal{R}$, we have 
 \[
   \|Y\|U_1U_1^{\top}\in\mathop{\arg\min}_{Z\in K\cap\mathcal{R}}\|Z-Y\|_F,
 \]
 which implies ${\rm dist}(Y,K\cap\mathcal{R})={\rm dist}(Y,\mathbb{S}_{+}^p\cap\mathcal{R})$. The conclusion then follows.
 \end{proof}
\begin{remark}\label{remark-noneig}
 The conclusion of Proposition \ref{dist-KpR} implies that for any $Y\in\mathbb{S}^p$, 
 \begin{align}\label{dist-ineq31}
 {\rm dist}(Y,K\cap\mathcal{R})&\le\|Y-\Pi_{\mathbb{S}^p\cap\mathbb{R}_{+}^{p\times p}}(Y)\|+{\rm dist}(\Pi_{\mathbb{S}^p\cap\mathbb{R}_{+}^{p\times p}}(Y),K\cap\mathcal{R})\nonumber\\
 &={\rm dist}(Y,\mathbb{S}^p\cap\mathbb{R}_{+}^{p\times p})+{\rm dist}(\Pi_{\mathbb{S}^p\cap\mathbb{R}_{+}^{p\times p}}(Y),\mathbb{S}_{+}^p\cap\mathcal{R})\nonumber\\
 &\le 2{\rm dist}(Y,\mathbb{S}^p\cap\mathbb{R}_{+}^{p\times p})+{\rm dist}(Y,\mathbb{S}_{+}^p\cap\mathcal{R}).
 \end{align}
 However, when replacing $\mathcal{R}$ with $\mathcal{R}_{r}:=\{Y\!\in\mathbb{S}^p\ |\ {\rm rank}(Y)\le r\}$ for $r\ge 2$, the inequality \eqref{dist-ineq31} does not necessarily hold, even locally. For example, consider $r=2$ and $\overline{Y}=\overline{v}^1(\overline{v}^1)^{\top}+\overline{v}^2(\overline{v}^2)^{\top}$ with $\overline{v}^1=(0,1,0,0)^{\top},\overline{v}^2=(0,1,0,1)^{\top}$. For each $k\in\mathbb{N}$, let $Y^k=v^{1,k}(v^{1,k})^{\top}+ v^{2,k}(v^{2,k})^{\top}$ with $v^{1,k}=(-\frac{1}{k},1,\frac{1}{k},0)^{\top}$ and $v^{2,k}=(\frac{1}{k},1,0,1)^{\top}$. Then, for each $k\in\mathbb{N}$,
 \[
  2{\rm dist}(Y^k,\mathbb{S}^p\cap\mathbb{R}_{+}^{p\times p})+{\rm dist}(Y^k,\mathcal{R}_{r}\cap\mathbb{S}_+^p)	={2\sqrt{2}}/{k^2}.
 \]
 Next we argue that there exists $c_1>0$ such that $\lim_{k\to\infty}k{\rm dist}(Y^k,K\cap\mathcal{R}_{r})\ge c_1$, so \eqref{dist-ineq31} does not hold when $\mathcal{R}$ is replaced with $\mathcal{R}_{r}$. Suppose on the contrary that $\lim_{k\to\infty}k{\rm dist}(Y^k,K\cap\mathcal{R}_{r})=0$. For each $k\in\mathbb{N}$, pick any $\overline{Y}^k\!\in\Pi_{K\cap\mathcal{R}_{r}}(Y^k)$. For each $k\in\mathbb{N}$, from \cite{Thomas74}, the nonnegative rank of $\overline{Y}^k$ is equal to its rank $2$, so there exist $u^{1,k}\in\mathbb{R}_{+}^4$ and $u^{2,k}\in\mathbb{R}_{+}^4$ such that $\overline{Y}^k\!=u^{1,k}(u^{1,k})^{\top}+u^{2,k}(u^{2,k})^{\top}$, and consequently,
 \begin{align}\label{temp-ineq32}
  \|Y^k-\overline{Y}^k\|_1
  &\ge (u_{1}^{1,k}u_{2}^{1,k}\!+\!u_{1}^{2,k}u_{2}^{2,k})+(u_{3}^{1,k}u_{4}^{1,k}\!+\!u_{3}^{2,k}u_{4}^{2,k})+\big|1\!-(u_{2}^{1,k}u_{4}^{1,k}\!+\!u_{2}^{2,k}u_{4}^{2,k})\big|\nonumber\\
  &\quad\ +\big|k^{-1}\!-(u_{2}^{1,k}u_{3}^{1,k}+u_{2}^{2,k}u_{3}^{2,k})\big|
  +\big|k^{-1}\!-(u_{1}^{1,k}u_{4}^{1,k}+u_{1}^{2,k}u_{4}^{2,k})\big|.
 \end{align}
 Obviously, $\|Y^k-\overline{Y}^k\|_1\ge |1-(u_{2}^{1,k}u_{4}^{1,k}\!+\!u_{2}^{2,k}u_{4}^{2,k})|$, which along with $\|Y^k-\overline{Y}^k\|_1\to 0$ as $k\to\infty$ implies that $(u_{2}^{1,k}u_{4}^{1,k}\!+\!u_{2}^{2,k}u_{4}^{2,k})\ge0.5$ for sufficiently large $k$. Without loss of generality, we assume $u_{2}^{2,k}u_{4}^{2,k}\ge0.25$ for large enough $k$. Then, there exists $c>0$ such that $\min\{u_{2}^{2,k},u_{4}^{2,k}\}\ge c$ for large enough $k$. Along with \eqref{temp-ineq32}, we have
 \begin{align*}
  \|Y^k-\overline{Y}^k\|_1
 &\ge (u_{1}^{1,k}u_{2}^{1,k}\!+\!u_{1}^{2,k}c)+(u_{3}^{1,k}u_{4}^{1,k}\!+\!u_{3}^{2,k}c)+|k^{-1}\!-(u_{2}^{1,k}u_{3}^{1,k}+u_{2}^{2,k}u_{3}^{2,k})|\\
 &\quad\ +|k^{-1}\!-(u_{1}^{1,k}u_{4}^{1,k}+u_{1}^{2,k}u_{4}^{2,k})|\quad{\rm for\ sufficiently\ large}\ k.
 \end{align*}
 Recalling that $0=\lim_{k\to{\infty}}k{\rm dist}(Y^k,K\cap\mathcal{R}_r)=\lim_{k\to{\infty}}k\|Y^k-\overline{Y}^k\|_1$, we have 
 \begin{subequations}
 \begin{align}\label{limit1}
 \lim_{k\to{\infty}}(u_{1}^{1,k}u_{2}^{1,k}\!+\!u_{1}^{2,k}c+u_{3}^{1,k}u_{4}^{1,k}\!+\!u_{3}^{2,k}c)k=0,\qquad\qquad\\
  \label{limit2}
 \lim_{k\to{\infty}}\big[|1-(u_{2}^{1,k}u_{3}^{1,k}+u_{2}^{2,k}u_{3}^{2,k})k|
 +|1-(u_{1}^{1,k}u_{4}^{1,k}+u_{1}^{2,k}u_{4}^{2,k})k|\big]=0.
 \end{align}
 \end{subequations}
 From \eqref{limit1} and the nonnegativity of $u^{1,k}$ and $u^{2,k}$, it follows $\lim_{k\to{\infty}}ku_{1}^{2,k}=0$ and $\lim_{k\to\infty}ku_{3}^{2,k}=0$, which along with \eqref{limit2} implies that $\lim_{k\to{\infty}}\big[|1-u_{2}^{1,k}u_{3}^{1,k}k|+|1-u_{1}^{1,k}u_{4}^{1,k}k|\big]=0$. Then, for sufficiently large $k$,  
 $u_{2}^{1,k}u_{3}^{1,k}\ge\frac{1}{2k}$ and $u_{1}^{1,k}u_{4}^{1,k}\ge\frac{1}{2k}$, so $k(u_{1}^{1,k}u_{2}^{1,k}+u_{3}^{1,k}u_{4}^{1,k})\ge\frac{u_{1}^{1,k}}{2u_{3}^{1,k}}+\frac{u_{3}^{1,k}}{2u_{1}^{1,k}}\ge 1$. On the other hand, the above \eqref{limit1} implies $\lim_{k\to{\infty}}[k(u_{1}^{1,k}u_{2}^{1,k}+u_{3}^{1,k}u_{4}^{1,k})]=0$. Thus, we get a contradiction. 
\end{remark}

Now we are in a position to establish the local error bounds stated in \eqref{aim-error}.
\begin{theorem}\label{Ebound-theorem}
 Let $\kappa:=(1+\!\sqrt{4n^2+n})(1\!+\!3\sqrt{2}n)$. For any $\overline{Y}\in\Gamma$, there exists $\varepsilon>0$ such that \eqref{aim-error} holds for all $Y\in\mathbb{B}(\overline{Y},\varepsilon)$, and consequently,     
 \begin{subequations}
 \begin{align}\label{distance-ineq1}
  {\rm dist}(Y,\Gamma)
  &\le 6\kappa\big[{\rm dist}(Y,T\cap\mathbb{L}\cap K)+{\rm dist}(Y,\mathcal{R})\big]\\
 \label{distance-ineq2}
  &\le6\kappa\big[{\rm dist}(Y,T\cap\Omega\cap K)+{\rm dist}(Y,\mathcal{R})\big]\\
 \label{distance-ineq3}
  &=6\kappa\big[{\rm dist}(Y,\Delta_0)+{\rm dist}(Y,\mathcal{R})\big].
 \end{align}
 \end{subequations}
\end{theorem}
\begin{proof}
 Let $\delta$ be the same as in Proposition \ref{dist-prop2}. Then, for any $Y\in\mathbb{B}(\overline{Y},\delta)$, 
 \begin{align}\label{temp-MF0}
  &(1\!+\!3\sqrt{2}n)^{-1}{\rm dist}(Y,\mathbb{L}\cap K\cap\mathcal{R}) \le {\rm dist}(Y,\mathbb{L})+{\rm dist}(Y,K\cap\mathcal{R})\nonumber\\
  &\stackrel{\eqref{dist-ineq31}}{\le}{\rm dist}(Y,\mathbb{L})+2{\rm dist}(Y,\mathbb{S}^p\cap\mathbb{R}_{+}^{p\times p})+{\rm  dist}(Y,\mathcal{R}\cap \mathbb{S}_{+}^p)\nonumber\\
  &\le{\rm dist}(Y,\mathbb{L})+2{\rm dist}(Y,\mathbb{S}^p\cap\mathbb{R}_{+}^{p\times p})+{\rm dist}(Y,\mathbb{S}_+^p)+{\rm dist}(\Pi_{\mathbb{S}_{+}^p}(Y),\mathcal{R})\nonumber\\
  &\le{\rm dist}(Y,\mathbb{L})+2{\rm dist}(Y,\mathbb{S}^p\cap\mathbb{R}_{+}^{p\times p})+2{\rm dist}(Y,\mathbb{S}_+^p)+{\rm dist}(Y,\mathcal{R})\nonumber\\
  &\le{\rm dist}(Y,\mathbb{L})+4{\rm dist}(Y,K)+{\rm dist}(Y,\mathcal{R}),
 \end{align}
 where the third inequality is obtained by using ${\rm dist}(Y,\mathcal{R}\cap\mathbb{S}_{+}^p)\le{\rm dist}(Y,\mathbb{S}_{+}^p)+{\rm dist}(\Pi_{\mathbb{S}_{+}^p}(Y),\mathcal{R}\cap\mathbb{S}_{+}^p)$ and ${\rm dist}(Z,\mathcal{R}\cap\mathbb{S}_{+}^p)={\rm dist}(Z,\mathcal{R})$ for all $Z\in\mathbb{S}^p_+$. Now set $\varepsilon$ to be the smaller $\delta$ from Propositions \ref{dist-prop1} and \ref{dist-prop2}. Pick any $Y\in\mathbb{B}(\overline{Y},\varepsilon)$. From Proposition \ref{dist-prop1} and the above \eqref{temp-MF0}, it follows that \eqref{aim-error} holds for this $\kappa$. Since $T\cap\mathbb{L}\cap K$ is contained in any of the sets $T,\mathbb{L}$ and $K$, the inequality \eqref{distance-ineq1} follows \eqref{aim-error}. The inequality \eqref{distance-ineq2} is due to the inclusion \eqref{imply1-set}, and the equality \eqref{distance-ineq3} follows $\Delta_0=T\cap\Omega\cap K$ by \eqref{imply1-set}. Thus, we complete the proof. 
\end{proof}
\begin{remark}\label{remark-ebound}
{\bf(a)} Theorem \ref{Ebound-theorem} provides the locally Lipschitzian error bounds for the set $\Gamma$ corresponding to its three expressions in \eqref{Eset-Gamma}, with the constant $\kappa$ depending only on the dimension $n$. According to \cite[Section 3.1]{Ioffe08}, these local error bounds actually imply a metric qualification condition for the set $\Gamma$ at any $\overline{Y}\in\Gamma$.

\noindent
{\bf(b)} In view of Remark \ref{remark-noneig}, for the relaxation set $\Delta_0\cap\mathcal{R}_{r}$ of $\Gamma$ with some $r\ge 2$, we cannot expect the locally Lipschitzian error bound as in \eqref{aim-error}. Now it is unclear whether locally H$\ddot{o}$lderian error bounds hold. We leave it for a future topic. 
\end{remark} 

By virtue of the compactness of $\Gamma$, we can establish its global error bounds as below. Unlike the above local error bounds, the involved constant $\kappa'$ is unavailable.
\begin{corollary}\label{global-ebound}
 For any bounded set $\Delta\subset\mathbb{S}^{p}$, there exists $\kappa'>0$ such that for all $Y\!\in\Delta$, the inequalities \eqref{aim-error} and \eqref{distance-ineq1}-\eqref{distance-ineq3} hold for $\kappa$ and $6\kappa$ replaced by $\kappa'$. 
\end{corollary}
\begin{proof}
 We present the detailed proof for the inequality \eqref{distance-ineq3}. Since the proofs for the other inequalities are similar, we omit them. By Theorem \ref{Ebound-theorem}, for each $\overline{Y}\in\Gamma\cap{\rm cl}(\Delta)$, where ${\rm cl}(\Delta)$ denotes the closure of $\Delta$, there exists $\delta_{\overline{Y}}>0$ such that for all $Y\in\mathbb{B}(\overline{Y},\delta_{\overline{Y}})$, the inequality \eqref{distance-ineq3} holds. Clearly, $\Gamma\cap{\rm cl}(\Delta)\subset\bigcup_{\overline{Y}\in\Gamma\cap{\rm cl}(\Delta)}\mathbb{B}^{\circ}(\overline{Y},\delta_{\overline{Y}})$ where $\mathbb{B}^{\circ}(\overline{Y},\delta_{\overline{Y}})$ is the interior of $\mathbb{B}(\overline{Y},\delta_{\overline{Y}})$. Since $\Gamma\cap{\rm cl}(\Delta)$ is compact, from the Heine-Borel covering theorem, there exist a finite number of points $\overline{Y}^1,\ldots,\overline{Y}^l\in\Gamma\cap{\rm cl}(\Delta)$ such that $\Gamma\cap{\rm cl}(\Delta)\subset\bigcup_{i=1}^l\mathbb{B}^{\circ}(\overline{Y}^i,\delta_{\overline{Y}^i}):=\aleph$. 
 
 Pick any $Y\in\Delta$. If $Y\in\Gamma\cup\aleph$, there exists $j\in[l]$ such that $Y\in\mathbb{B}^{\circ}(\overline{Y}^j,\delta_{\overline{Y}^j})$, so the point $Y$ satisfies \eqref{distance-ineq3}. Next we consider that $Y\notin\Gamma\cup\aleph$. Let $\widetilde{\Delta}={\rm cl}[\Delta\backslash(\Gamma\cup\aleph)]$. Clearly, $Y\in\widetilde{\Delta}\subset{\rm cl}(\Delta)$. Furthermore, $\widetilde{\Delta}\cap\aleph=\emptyset$ (if not, there exists $\widehat{Z}\in\widetilde{\Delta}\cap \aleph$. By the definition of $\widetilde{\Delta}$, there is a sequence $\{\widehat{Z}^k\}_{k\in\mathbb{N}}\subset\Delta\backslash(\Gamma\cup\aleph)$ such that $\lim_{k\to\infty}\widehat{Z}^k=\widehat{Z}$. From the openness of $\aleph$ and $\widehat{Z}\in\aleph$, we have $\widehat{Z}^k\in\aleph$ for sufficiently large $k$, a contradiction to $\widehat{Z}^k\notin\aleph$ for all $k$). Now we claim that there exists $\mu>0$ such that 
\begin{equation}\label{dist-mu}
  \min_{Z\in\widetilde{\Delta}}\,{\rm dist}(Z,\Delta_0)+{\rm dist}(Z,\mathcal{R})\ge \mu.
\end{equation}
If not, we can find a sequence $\{Z^k\}_{k\in\mathbb{N}}\subset\widetilde{\Delta}$ such that ${\rm dist}(Z^k,\Delta_0)+{\rm dist}(Z^k,\mathcal{R})\le\frac{1}{k}$ for each $k$. By the compactness of $\widetilde{\Delta}$ and the continuity of the distance function, the sequence $\{Z^k\}_{k\in\mathbb{N}}$ has a cluster point $\overline{Z}\in\widetilde{\Delta}$ satisfying ${\rm dist}(\overline{Z},\Delta_0)+{\rm dist}(\overline{Z},\mathcal{R})\le 0$. The latter implies  $\overline{Z}\in\Delta_0\cap\mathcal{R}=\Gamma$. Along with  $\overline{Z}\in\widetilde{\Delta}\subset{\rm cl}(\Delta)$, it follows $\overline{Z}\in\Gamma\cap{\rm cl}(\Delta)\subset\aleph$. Recall that $\widetilde{\Delta}\cap\aleph=\emptyset$. Then, $\overline{Z}\notin\widetilde{\Delta}$, a contradiction to $\overline{Z}\in\widetilde{\Delta}$. The claimed \eqref{dist-mu} holds. Recall that $Y\in\widetilde{\Delta}$, so ${\rm dist}(Y,\Delta_0)+{\rm dist}(Y,\mathcal{R})\ge\mu$. 
In addition, from the compactness of $\widetilde{\Delta}$ and $\Gamma$, there exists a constant $c>0$ such that ${\rm dist}(Z,\Gamma)\le c$ for all $Z\in\widetilde{\Delta}$. Then, ${\rm dist}(Y,\Gamma)\le c=(c/\mu)\mu\le (c/\mu)[{\rm dist}(Y,\Delta_0)+{\rm dist}(Y,\mathcal{R})]$. Thus, \eqref{distance-ineq3} holds for all $Y\in\Delta$ with $6\kappa$ replaced by $\kappa'=\max\{6\kappa,c/\mu\}$.    
\end{proof}
\section{Global exact penalties}\label{sec4}

In this section, we apply the error bounds to prove that the penalty problem \eqref{gDNN-penalty} is a global exact penalty of \eqref{gDNN-mpec}, and that the BM factorization \eqref{gfac-penalty} is a global exact penalty of \eqref{gfac-mpec}. First, the following proposition states that the multifuntions $\Upsilon,\Upsilon_1$ and $\Upsilon_0$ are locally upper Lipschitzian at $\tau=0$. 
\begin{proposition}\label{calmness-prop}
The multifunctions $\Upsilon,\Upsilon_{\!1}$ and $\Upsilon_0$  are locally upper Lipschitzian at $0$ with modulus $\kappa'$, and $\mathbb{R}\ni\tau\rightrightarrows\widehat{\Upsilon}(\tau):=\big\{Y\in\Sigma\ |\ {\rm tr}(Y)-\interleave Y\interleave=\tau\big\}$ with $\interleave \cdot\interleave=\|\cdot\|_F$ is locally upper Lipschitzian at $0$ with modulus $2\kappa'$.
\end{proposition}
\begin{proof}	
 We first prove the conclusion for the mapping $\Upsilon$.   
 Fix any $\delta>0$. Consider any $\tau\in\mathbb{R}$ with $|\tau|\le\delta$. We claim that $\Upsilon(\tau)\subset\Upsilon(0)+\kappa'\tau\mathbb{B}$ holds. Indeed, when $\tau< 0$, since $\Upsilon(\tau)=\emptyset$, the desired inclusion is trivial, so it suffices to consider that $0\le\tau\le\delta$ with $\Upsilon(\tau)\ne\emptyset$. Pick any $Y\in\Upsilon(\tau)$. From the definition $\Upsilon$, $Y\in T\cap\mathbb{L}\cap K$ and ${\rm tr}(Y)-\interleave Y\interleave=\tau$. Invoking Corollary \ref{global-ebound} with $\Delta=T\cap\mathbb{L}\cap K$ leads to 
 \begin{align}\label{temp-ULineq}
 {\rm dist}(Y,\Upsilon(0))&={\rm dist}(Y,\Gamma)\le \kappa'\big[{\rm dist}(Y,T\cap\mathbb{L}\cap K)+{\rm dist}(Y,\mathcal{R})]\nonumber\\
 &=\kappa'{\rm dist}(Y,\mathcal{R})\le\kappa'\big[{\rm tr}(Y)-\| Y\|\big]=\kappa'\tau.
 \end{align}
 This implies the desired inclusion, so $\Upsilon$ is locally upper Lipschitzian at $0$. Using the same arguments can prove the local upper Lipschitz property of $\Upsilon_1$ and $\Upsilon_0$ at $0$. 

 Next we prove that the mapping $\widehat{\Upsilon}$ with $\Sigma=T\cap\mathbb{L}\cap K$ and $\interleave\cdot\interleave=\|\cdot\|_F$ is locally upper Lipschitzian at $0$. Consider any $\tau\in\mathbb{R}$ with $|\tau|\le\delta$. If $\widehat{\Upsilon}(\tau)=\emptyset$, then 
 \begin{equation}\label{temp-inclusion}
  \widehat{\Upsilon}(\tau)\subset\widehat{\Upsilon}(0)+2\kappa'\tau\mathbb{B} 
\end{equation}
automatically hold. Next we consider the case that $\widehat{\Upsilon}(\tau)\ne\emptyset$. Pick any $Y\in\widehat{\Upsilon}(\tau)$. Obviously, $Y\in T\cap\mathbb{L}\cap K$. Then, $Y\in\Upsilon(\omega)$ with $\omega={\rm tr}(Y)-\|Y\|_F\ge 0$. Consequently, 
 \begin{align*}
 {\rm dist}(Y,\widehat{\Upsilon}(0))&={\rm dist}(Y,\Gamma)\stackrel{\eqref{temp-ULineq}}{\le} \kappa'\big[{\rm tr}(Y)-\| Y\|\big]\le 2\kappa'\big[{\rm tr}(Y)-\| Y\|_F\big]=2\kappa'\tau,
 \end{align*}
 where the second inequality is due to \cite[Lemma 2]{QianPanLiu23} with $r=1$. This shows that the inclusion \eqref{temp-inclusion} holds for this case. We conclude that  $\widehat{\Upsilon}$ with $\Sigma=T\cap\mathbb{L}\cap K$ is locally upper Lipschitzian at $0$ with modulus $2\kappa'$. Similarly, when $\Sigma=T\cap\Omega\cap K$ and $\Delta_0$, the conclusion also holds for the mapping $\widehat{\Upsilon}$. Thus, we complete the proof.
 \end{proof}
 
 Notice that Proposition \ref{calmness-prop} implies the calmness of the mappings $\Upsilon,\Upsilon_1$ and $\Upsilon_0$ at $\tau=0$ for all $Y\in\Gamma$. Then, invoking Lemma \ref{lemma-calm} results in the partial calmness of the problem \eqref{gDNN-mpec} on the set of local optimal solutions, so that of the problem \eqref{gDNN-mpec} on the set of global optimal solutions. The latter along with \cite[Proposition 2.1(b)]{LiuBiPan18} and the compactness of $\Sigma$ implies the global exact penalty result of \eqref{gDNN-penalty}. 
\begin{theorem}\label{penalty-theorem}
 The problem \eqref{gDNN-mpec} is partially calm on the set of local optimal solutions when a perturbation is imposed on the constraint ${\rm tr}(Y)-\interleave Y\interleave=0$, and there exists $\overline{\rho}>0$ such that the penalty problem \eqref{gDNN-penalty} associated with every $\rho>\overline{\rho}$ 
 has the same global optimal solution set as \eqref{gDNN-mpec}. Consequently, the penalty problem \eqref{EroneDNNcon-penalty} is a global exact penalty of \eqref{EroneDNNcon}, so is the following penalty problem of \eqref{EroneDNNcon1}:
 \begin{align*}
 &\min_{Y\in\mathbb{S}^p}\ \langle G,Y\rangle+\rho({\rm tr}(Y)-\|Y\|)\nonumber\\
 &\ {\rm s.t.}\ \ \langle I,Y^{ii}\rangle=1,\, {\textstyle\sum_{t=1}^n}{(Y^{tt})_{ii}}=1\ \ {\rm for}\  i\in[n],\nonumber\\
 &\qquad\ \langle D,Y\rangle=0,\ \langle ee^{\top},Y\rangle=p,\,Y\!\in K.
 \end{align*}
 \end{theorem}

 Next we establish the partial calmness of \eqref{gfac-mpec} on its local optimal solution set, and get its global exact penalty induced by the DC constraint $\|V\|_F^2-\interleave V^{\top}V\interleave =0$.
\begin{theorem}\label{pcalmV1}
 The problem \eqref{gfac-mpec} is partially calm on the set of local optimal solutions when a perturbation is imposed on the constraint $\|V\|_F^2-\interleave V^{\top}V\interleave=0$, and there exists $\overline{\rho}>0$ such that the penalty problem \eqref{gfac-penalty} associated with every $\rho\ge\overline{\rho}$ has the same global optimal solution set as \eqref{gfac-mpec}. Consequently, the penalty problem
 \begin{align}\label{Epfactor-ding}
  &\min_{V\in\mathbb{R}^{m\times p}}\,\langle G,V^{\top}V\rangle+\rho(\|V\|_F^2-\|V^{\top}V\|)\nonumber\\
  &\quad\ {\rm s.t.}\quad\|V_i\|_F^2=1,\,{\textstyle\sum_{t=1}^n}\|(V_t)_i\|^2=1\ \ \forall i\in[n],\\
  &\qquad\quad\ \ \langle D,V^{\top}V\rangle=0,\,\langle ee^{\top},\,V^{\top}V\rangle=p,\ V^{\top}V\in\mathbb{R}_{+}^{p\times p}\nonumber
 \end{align}
 is a global exact penalty for the BM factorization of \eqref{EroneDNNcon1}, which takes the form of
 \begin{equation*}
  \min_{V\in\mathbb{R}^{m\times p}}\big\{\langle G,V^{\top}V\rangle\ \ {\rm s.t.}\ \ \|V\|_F^2-\|V^{\top}V\|=0,\,V^{\top}V\in T\cap\Omega\cap\mathbb{R}_{+}^{p\times p}\big\}.
 \end{equation*}
 \end{theorem}
 \begin{proof}
  Define a partial perturbation mapping for the feasible set of  \eqref{gfac-mpec} as
 \begin{equation}\label{Xi-map}
  \mathcal{E}(\tau):=\Big\{V\in\mathbb{R}^{m\times p}\ |\
  \|V\|_F^2-\interleave V^{\top}V\interleave=\tau,\,V^{\top}V\in\Xi\Big\}.
 \end{equation}
 Denote by $\mathcal{V}^*$ the set of local optimal solutions of \eqref{gfac-mpec}. Pick any $V^*\!\in\mathcal{V}^*$. Then, $Y^*=({V^*})^{\top}V^*\in\Gamma$. By virtue of Lemma \ref{Gamma-set}, $\Gamma$ is a discrete set, so $\Gamma$ coincides with the set of local optimal solutions of \eqref{gDNN-mpec}. From Theorem \ref{penalty-theorem}, the problem \eqref{gDNN-mpec} is partially calm on $\Gamma$. Then there exist $\varepsilon'\in(0,1)$ and $\mu>0$ such that for all $\tau\in\mathbb{R}$ and all $Y\in\mathbb{B}(Y^*,\varepsilon')\cap\widehat{\Upsilon}(\tau)$, where $\widehat{\Upsilon}$ is the mapping appearing in Proposition \ref{calmness-prop}, 
 \begin{equation}\label{fineq1}
  f(Y)-f(Y^*)+\mu\big({\rm tr}(Y)-\interleave Y\interleave\big)\ge0.
 \end{equation}
 Let $\varepsilon=\frac{\varepsilon'}{2\|V^*\|+1}$. Pick any $\tau\in\mathbb{R}$ and  $V\in\mathbb{B}(V^*,\varepsilon)\cap\mathcal{E}(\tau)$. Then, it holds 
 \[
   \|V^{\top}V-(V^*)^{\top}V^*\|_F\le(\|V\|+\|V^*\|)\|V-V^*\|_F\le(2\|V^*\|+\varepsilon)\varepsilon\le\varepsilon'.
 \]
 This shows that $V^{\top}V\in\mathbb{B}(Y^*,\varepsilon')$. In addition, combining $V\in\mathcal{E}(\tau)$ with the definitions of $\mathcal{E}$ and $\widehat{\Upsilon}$, we have $V^{\top}V\in\widehat{\Upsilon}(\tau)$. The two sides show that $V^{\top}V\in\mathbb{B}(Y^*,\varepsilon')\cap\widehat{\Upsilon}(\tau)$.  Consequently, the inequality \eqref{fineq1} holds for $Y=V^{\top}V$, that is,
 \[
  f(V^{\top}V)-f((V^*)^{\top}V^*)+\mu(\|V\|_F^2-\interleave V\interleave)\ge 0.
 \]
 The conclusion follows the arbitrariness of $\tau\in\mathbb{R}$ and $V\in\mathbb{B}(V^*,\varepsilon)\cap\mathcal{E}(\tau)$. 
 \end{proof}

 It is worth pointing out that the calmness of $\mathcal{E}$ at $0$ for any $\overline{V}\in\mathcal{E}(0)$ is generally not true though the calmness of $\widehat{\Upsilon}$ at $0$ for $\overline{V}^{\top}\overline{V}$ always holds. For example, consider 
 \[
 \overline{V}=\left(\begin{matrix}
	1 & 0 & 0 &  1 \\
	0 & 0 & 0 & 0
 \end{matrix}\right)\ {\rm and}\ V_{\tau}=\left(\begin{matrix}
	\sqrt{1\!-\!\tau/2} & 0 & 0 &  \sqrt{1\!-\!\tau/2} \\
	0 & \sqrt{\tau/2} & \sqrt{\tau/2} & 0
\end{matrix}\right)\ \ {\rm for}\ \tau\in(0,2). 
\]
Consider the mapping $\mathcal{E}$ in \eqref{Xi-map} with $\interleave \cdot\interleave=\|\cdot\|$. 
One can check that $\overline{V}\in\mathcal{E}(0)$, since $\|\overline{V}\|_{F}^2=\|\overline{V}^{\top}\overline{V}\|$. For any $\tau\in(0,2)$, $V_{\tau}\in\mathcal{E}(\tau)$ because $\|V_{\tau}\|_F^2=2$ and $\|V_{\tau}^{\top}V_{\tau}\|=2(1\!-\!\tau/2)$. However, 
${\rm dist}(V_{\tau},\mathcal{E}(0))\!=\!\|V_{\tau}\!-\overline{V}\|_F=O(\sqrt{\tau})$. 
 Clearly, there is no constant $\kappa>0$ such that $\mathcal{E}(\tau)\cap\mathbb{B}(\overline{V},1/2)\subset\mathcal{E}(0)+\kappa\tau\mathbb{B}$.

\section{Relaxation approach based on global exact penalty}\label{sec5}

In the last section, the penalty problem \eqref{gDNN-penalty} associated with every $\rho>\overline{\rho}$ is proved to possess the same global optimal solution set as \eqref{gDNN-mpec}. The same conclusion also holds for their BM factorizations \eqref{gfac-penalty} and \eqref{gfac-mpec}. This section illustrates their application in developing relaxation approaches to seek a rank-one approximate feasible solution of the QAP, by using the penalty problem \eqref{gfac-penalty} with $f(\cdot)=\langle G,\cdot\rangle,\Xi=T\cap\Omega\cap\mathbb{R}_{+}^{p\times p}$ and $\interleave\cdot\interleave=\|\cdot\|_F$. This penalty problem, according to \eqref{Tset-def}-\eqref{Omega-def}, has the following form  
\begin{align}\label{gfac-penalty2}
 &\min_{V\in\mathbb{R}^{m\times p}}\langle G,V^{\top}V\rangle+\rho(\|V\|_F^2-\|V^{\top}V\|_F)\nonumber\\
 &\quad\ {\rm s.t.}\ \ V^{\top}V\in\mathcal{Y},\ \langle E,V^{\top}V\rangle=1\ \ {\rm with}\ E:=p^{-1}ee^{\top},
\end{align}
where $\mathcal{Y}:=\big\{Y\!\in\mathbb{S}^{p}\cap\mathbb{R}_{+}^{p\times p}\,|\,\langle D,Y\rangle=0,{\textstyle\sum_{t=1}^n}{(Y^{tt})_{ii}}=1,\langle I,Y^{ii}\rangle=1\ \ \forall i\in[n]\big\}$. Compared with the previous \eqref{Epfactor-ding}, the problem \eqref{gfac-penalty2} has an advantage in the smoothness of its objective function on the feasible set. Since the threshold $\overline{\rho}$ is unknown and every penalty problem is nonconvex, one cannot expect a high-quality solution when solving a single penalty problem \eqref{gfac-penalty2} with numerical algorithms. Inspired by this, we propose a relaxation approach by computing approximate stationary points of a finite number of penalty problems with increasing penalty factor $\rho$. 
 \begin{algorithm}[H]
 \caption{\label{Alg-factor}{\bf(Relaxation approach based on \eqref{gfac-penalty2})}}
 \textbf{Input:} $m\in[p],\epsilon_1\in(0,1),\epsilon_2\!\in(0,1),\underline{\tau}\in(0,1),\tau_1\in(\underline{\tau},1),\varsigma\in(0,1),\rho_{\rm max}>0$, \hspace*{1.3cm} $0<\rho_1<\rho_{\rm max},\sigma>1$ and $V^0\in\mathbb{R}^{m\times p}$. 
 
 \noindent
 \textbf{For} $l=1,2,\ldots$
 \begin{itemize}
 \item [1.] Starting from $V^{l-1}$, seek a $\tau_{l}$-approximate stationary point $V^{l}$ of the penalty problem \eqref{gfac-penalty2} associated with $\rho=\rho_{l}$. 
		
\item [2.] If $\|V^{l}\|_F^2-\|(V^{l})^{\top}V^{l}\|_{F}\le\epsilon_1$ and $\sqrt{[{\rm dist}((V^{l})^{\top}V^{l},\mathcal{Y})]^2\!+\!|\langle E,(V^{l})^{\top}V^{l}\rangle\!-\!1|^2}\le \epsilon_2$, stop; else go to step 3. 
		
 \item[3.] Let $\rho_{l+1}\leftarrow\min\{\sigma\rho_{l},\rho_{\rm max}\}$ and $
 \tau_{l+1}\leftarrow\max\{\varsigma\tau_{l},\underline{\tau}\}$.
 \end{itemize}
 \textbf{end (For)}
 \end{algorithm}
\begin{remark}\label{remark-Alg}
 A matrix $V\in\mathbb{R}^{m\times p}$ is called a $\tau_{l}$-approximate stationary point of the problem \eqref{gfac-penalty2} if there exist $Y\in\mathcal{Y}$ and $(S,\lambda)\in\mathcal{N}_{\mathcal{Y}}(Y)\times\mathbb{R}$ such that 
 \begin{subequations}
\begin{align}\label{cpoint-equa1}
 & \sqrt{\|V^{\top}V-Y\|_F^2+|\langle E,V^{\top}V\rangle-1|^2}\le\tau_{l},\\
 \label{cpoint-equa2}
 &\big\|2V\big[G+\rho I-\rho\|V^{\top}V\|_F^{-1}V^{\top}V+\lambda E+S\big]\big\|_F\le 10\tau_{l},
 \end{align}         
\end{subequations}
 where $\mathcal{N}_{\mathcal{Y}}(Y)$ is the normal cone to the closed convex set $\mathcal{Y}$ at $Y$. We see that, if $V^{l}$ is a rank-one $\tau_{l}$-approximate stationary point of \eqref{gfac-penalty2} with $\rho=\rho_{l}$ and $\tau_{l}\le\epsilon_2$, Algorithm \ref{Alg-factor} must stop at the $l$th iteration. Now one can use $V^{l}$ to construct a rank-one DNN matrix satisfying approximately the constraints of \eqref{EroneDNNcon1} by Proposition \ref{prop4.1}.
\end{remark}
 \begin{proposition}\label{prop4.1}
 Let $V^{l_f}$ be the output of Algorithm \ref{Alg-factor}, and define the rank-one DNN matrix $\widetilde{Y}^{l_f}\!:=\|V^{l_f}\|^2|Q_1||Q_1|^{\top}$, where $Q_1$ is the first column of any orthogonal matrix $Q\in\mathbb{O}^{p}((V^{l_f})^{\top}V^{l_f})$. If ${\rm rank}(V^{l_f})=1$ and $0<\tau_{l_f}<\frac{n}{n+m}$, then it holds 
 \[
   {\rm dist}(\widetilde{Y}^{l_f},\mathcal{Y})\le \Big[1+\frac{8m\|V^{l_{f}}\|^2}{n-(n+m)\tau_{l_{f}}}\Big]\tau_{l_{f}}.
 \]
\end{proposition}
\begin{proof}
 Let $Y^{l_f}\!:=(V^{l_f})^{\top}V^{l_f}$. By Remark \ref{remark-Alg}, ${\rm dist}(Y^{l_{f}},\mathcal{Y})\le\tau_{l_{f}}$. Recall that $\mathcal{Y}\subset T$. The matrix $\overline{Y}^{l_{f}}\!:=\Pi_{T}(Y^{l_{f}})$ satisfies $\|\overline{Y}^{l_{f}}\!-Y^{l_{f}}\|_F\le {\rm dist}(Y^{l_f},\mathcal{Y})\le\tau_{l_f}$. Then  
 \begin{align*}
  |\|V^{l_{f}}\|_F^2-n|&=|\langle I,Y^{l_f}\rangle-n|
  =|\langle I,Y^{l_f}-\overline{Y}^{l_{f}}\rangle+\langle I,\overline{Y}^{l_{f}}\rangle-n|\\
 &\le n\|\overline{Y}^{l_{f}}\!-\!Y^{l_{f}}\|_F+|\langle I,\overline{Y}^{l_f}\rangle-n|=n\|\overline{Y}^{l_{f}}\!-\!Y^{l_{f}}\|_F\le n\tau_{l_{f}}, 
 \end{align*}
 where the third equality is due to $\overline{Y}^{l_{f}}\in T$ and the expression of $T$. Together with $\|V^{l_{f}}\|_F^2\le m\|V^{l_{f}}\|^2$, it follows that $\|V^{l_{f}}\|^2-\tau_{l_{f}}\ge\frac{1}{m}[n-(n+m)\tau_{l_f}]>0$. 
 
 Let $Y^{+}\!:=\Pi_{\mathbb{R}_{+}^{p\times p}}(Y^{l_f})$. From $\mathcal{Y}\subset\mathbb{R}_{+}^{p\times p}$, we get $\|Y^{l_{f}}-Y^{+}\|_F={\rm dist}(Y^{l_{f}},\mathbb{R}_{+}^{p\times p})\le{\rm dist}(Y^{l_{f}},\mathcal{Y})\le \tau_{l_f}$, which implies $\|Y^{l_{f}}-Y^{+}\|\le\tau_{l_f}$. Consequently,  
 \begin{equation}\label{temp-Vlf}
  0<\frac{1}{m}[n-(n+m)\tau_{l_f}]\le\|V^{l_{f}}\|^2-\tau_{l_f}\le\|Y^{+}\|\le \|V^{l_{f}}\|^2+\tau_{l_f}. 
  \end{equation}
  By the proof of Proposition \ref{dist-KpR}, there exists $U_1\in\mathbb{R}_{+}^{p}\backslash\{0\}$ such that $U_1$ is an eigenvector of $Y^{+}$ associated with $\|Y^{+}\|$. From \cite[Theorem 2.1]{Dopico00} with $A=Y^{l_f},\widetilde{A}=Y^{+}$, 
 \[
  \min_{\alpha\in\{-1,1\}}\|Q_1\alpha-U_1\|\le \frac{2\tau_{l_{f}}}{\|Y^{+}\|}\ \le\frac{2\tau_{l_{f}}}{\|V^{l_{f}}\|^2-\tau_{l_f}}. 
 \] 
 Let $\alpha^*$ be an optimal solution of the above minimization. Along with $U_1\in\mathbb{R}_{+}^{p}\backslash\{0\}$,  
 \begin{equation}\label{perturb-ineq1}
   \||Q_1|-U_1\|\le\|Q_1\alpha^*-U_1\|\le \frac{2\tau_{l_{f}}}{\|V^{l_{f}}\|^2-\tau_{l_f}}.
 \end{equation}
 Since ${\rm rank}(V^{l_f})=1$, it follows $Y^{l_f}=\|V^{l_f}\|^2Q_1Q_1^{\top}=\|V^{l_f}\|^2(Q_1\alpha^*)(Q_1\alpha^*)^{\top}$. Then,
 \begin{align*}
 \|Y^{l_{f}}-\widetilde{Y}^{l_f}\|_F&=\|V^{l_f}\|^2\|(Q_1\alpha^*-|Q_1|)(Q_1\alpha^*)^{\top}-|Q_1|(|Q_1|-Q_1\alpha^*)^{\top}\|_F\\
 &\le 2\|V^{l_f}\|^2\|Q_1\alpha^*-|Q_1|\|\\
 &\le2\|V^{l_f}\|^2\|Q_1\alpha^*-U_1\|+2\|V^{l_f}\|^2\|U_1-|Q_1|\|\\
 &\stackrel{\eqref{perturb-ineq1}}{\le} \frac{8\tau_{l_{f}}\|V^{l_{f}}\|^2}{\|V^{l_{f}}\|^2-\tau_{l_f}}\stackrel{\eqref{temp-Vlf}}{\le} \frac{8m\|V^{l_{f}}\|^2\tau_{l_{f}}}{n-(n+m)\tau_{l_{f}}}.
 \end{align*}
 This, together with ${\rm dist}(Y^{l_{f}},\mathcal{Y})\le\tau_{l_{f}}$, implies the desired conclusion.  
 \end{proof} 

\subsection{ALM for solving penalty subproblems}\label{sec5.1}

We propose an ALM for seeking a $\tau_{l}$-approximate stationary point of the problem \eqref{gfac-penalty2} with $\rho=\rho_{l}$. With an additional variable $Y\in\mathbb{S}^p$, the problem \eqref{gfac-penalty2} with $\rho=\rho_{l}$ is equivalently written as 
\begin{align}\label{Efac-penalty22}
 &\min_{V\in\mathbb{R}^{m\times p},Y\in\mathcal{Y}}\langle G,V^{\top}V\rangle+\rho_{l}(\|V\|_F^2-\|V^{\top}V\|_F)\nonumber\\
 &\qquad {\rm s.t.}\ \ Y-V^{\top}V=0,\ \langle E,V^{\top}V\rangle=1.
\end{align}
 For a given $\beta>0$, the augmented Lagrangian function of \eqref{Efac-penalty22} takes the form of  
 \begin{align*}
 L_{\rho_{l},\beta}(V,Y;S,\lambda)&:=\langle G,V^{\top}V\rangle+\rho_{l}(\|V\|_F^2\!-\!\|V^{\top}V\|_F)+\chi_{\mathcal{Y}}(Y)+\lambda\big(\langle E,V^{\top}V\rangle\!-\!1\big)\\
 &\quad\ +\langle S,V^{\top}V\!-\!Y\rangle+\frac{\beta}{2}\big[(\langle E,V^{\top}V\rangle-1)^2+\|Y-V^{\top}V\|_F^2\big],
 \end{align*} 
 where $S\in\mathbb{S}^p$ and $\lambda\in\mathbb{R}$ are the Lagrange multiplier corresponding to the constraints $Y=V^{\top}V$ and $\langle E,V^{\top}V\rangle=1$, respectively. 
 For any given $\beta>0$ and $(S,\lambda)\in\mathbb{S}^p\times\mathbb{R}$, the basic iteration step of the ALM is 
 $(V^{+},Y^{+})\approx \mathop{\arg\min}_{V,Y}L_{\rho_{l},\beta}(V,Y;S,\lambda)$. Let $\mathcal{L}_{\rho_l,\beta}(V,S,\lambda):=\min_{Y\in\mathbb{S}^{p}}L_{\rho_{l},\beta}(V,Y;S,\lambda)$. An elementary calculation results in 
 \begin{align*}
 \mathcal{L}_{\rho_l,\beta}(V,S,\lambda)&=\langle G,V^{\top}V\rangle+\!\rho_{l}(\|V\|_F^2\!-\!\|V^{\top}V\|_F)
   +\!\lambda(\langle E,V^{\top}V\rangle\!-\!1)-\frac{1}{2\beta}\|S\|_F^2\\
 &\quad\ \ +\!\frac{\beta}{2}(\langle E,V^{\top}V\rangle\!-\!1)^2+\frac{\beta}{2}\big\|V^{\top}V+\beta^{-1}S-\Pi_{\mathcal{Y}}(V^{\top}V+\beta^{-1}S)\big\|_F^2
 	.
\end{align*} 
Then, the basic iteration step of the ALM is simplified as follows 
 \[
   V^{+}\approx \mathop{\arg\min}_{V\in\mathbb{R}^{m\times p}}\mathcal{L}_{\rho_l,\beta}(V,S,\lambda)\ \ {\rm and}\ \ Y^{+}=\Pi_{\mathcal{Y}}((V^{+})^{\top}V^{+}+\beta^{-1}S).
 \]
 At the $k$th iteration, our ALM seeks a $\varepsilon_k$-critical point of $\mathcal{L}_{\rho_l,\beta_{k}}(\cdot,\widehat{S}^{l,k},\widehat{\lambda}^{l,k})$ with 
 \[
  (\widehat{S}^{l,k},\widehat{\lambda}^{l,k})\in\mathbb{B}_{\varpi}:=\Big\{(S,\lambda)\in\mathbb{S}^p\times\mathbb{R}\,|\, \sqrt{\|S\|_F^2+\lambda^2}\le\varpi\Big\}\ \ {\rm for\ some}\ \varpi>0,
 \]
 and then update the multipliers. As will be shown in Proposition \ref{ALM-feasible}, the modified multiplier pair $(\widehat{S}^{l,k},\widehat{\lambda}^{l,k})$ is the key to achieve the asymptotic feasibility of the ALM.
 \renewcommand{\thealgorithm}{A}
 \begin{algorithm}[H]
 \caption{\label{ALM}({\bf ALM for the penalty subproblems})}
 \textbf{Input:} $\rho_{l}>0,\varepsilon=10\tau_l,\overline{\beta}>0,\gamma>1,\beta_0>0$ and $(S^{l,0},\lambda^{l,0})$. Let $V^{l,0}\!=V^{l-1}$ and \hspace*{1.2cm} $\varpi=\max\{10^3,|\langle G,(V^{l,0})^{\top}V^{l,0}\rangle|\}$.

 \noindent
 \textbf{For} $k=0,1,2,\ldots$
 \begin{itemize}
  \item[1.] Compute the projection  $(\widehat{S}^{l,k},\widehat{\lambda}^{l,k})$ of $(S^{l,k},\lambda^{l,k})$ onto the closed ball $\mathbb{B}_{\varpi}$.
		
 \item[2.] Seek an $\varepsilon$-critical point $V^{l,k+1}$ of the function $\mathcal{L}_{\rho_l,\beta_k}(\cdot,\widehat{S}^{l,k},\widehat{\lambda}^{l,k})$, and compute 
 \[
  Y^{l,k+1}=\Pi_{\mathcal{Y}}((V^{l,k+1})^{\top}V^{l,k+1}+\beta_k^{-1}\widehat{S}^{l,k}).
 \]
		
\item[3.] Update the Lagrange multipliers by the following formulas
\begin{subequations}
\begin{align}\label{multiplier-S} 
 S^{l,k+1}&=\widehat{S}^{l,k}+\beta_k((V^{l,k+1})^{\top}V^{l,k+1}-Y^{l,k+1}),\\
 \label{multiplier-lambda} 
 \lambda^{l,k+1}&=\widehat{\lambda}^{l,k}+\beta_k(\langle E,(V^{l,k+1})^{\top}V^{l,k+1}\rangle-1).
\end{align}
\end{subequations}
 		
\item[4.] Set $\beta_{k+1}:=\min\{\overline{\beta},\gamma\beta_k\}$.
\end{itemize}
\textbf{end (For)}
\end{algorithm}
\begin{remark}\label{remark-ALM}
 {\bf(a)} Recall that the matrices $A,B$ and $C$ are assumed to be nonnegative, so is the matrix $G$. From the expression of $\mathcal{L}_{\rho_l,\beta_k}(\cdot,\widehat{S}^{l,k},\widehat{\lambda}^{l,k})$, we infer that it is coercive, which means that $\min_{V\in\mathbb{R}^{m\times p}}\mathcal{L}_{\rho_l,\beta_k}(V,\widehat{S}^{l,k},\widehat{\lambda}^{l,k})$ has an optimal solution. Consequently, $V^{l,k+1}$ in step 2 exists and Algorithm \ref{ALM} is well defined. 
	
\noindent
{\bf(b)} By the expression of the set $\mathcal{Y}$, for any $Y\in\mathbb{S}^p$, the matrix ${\rm mat}({\rm diag}(\Pi_{\mathcal{Y}}(Y)))$ is precisely the projection of ${\rm mat}({\rm diag}(Y))$ onto $\mathcal{D}$, the convex hull of $\mathcal{P}$, so we can apply the fast semismooth Newton method developed in \cite{LiSun20} to compute the diagonal elements of $\Pi_{\mathcal{Y}}(Y)$. The non-diagonal elements of $\Pi_{\mathcal{Y}}(Y)$ have a closed form. Thus, the computation of the projection in step 3 can be finished with the solver in \cite{LiSun20}. 
 
\noindent
{\bf(c)} From the expression of $\mathcal{L}_{\rho_l,\beta_{k-1}}(\cdot,\widehat{S}^{l,k-1},\widehat{\lambda}^{l,k-1})$ and $\varepsilon=10\tau_{l}$, when $V^{l,k}$ is a $\varepsilon$-critical point of the function $\mathcal{L}_{\rho_l,\beta_{k-1}}(\cdot,\widehat{S}^{l,k-1},\widehat{\lambda}^{l,k-1})$, it holds that  
\begin{equation*}
 \big\|2V^{l,k}\big[G+\rho_{l}I-\rho_{l}\|(V^{l,k})^{\top}V^{l,k}\|_F^{-1}(V^{l,k})^{\top}V^{l,k}+\lambda^{l,k}E+S^{l,k}\big]\big\|_F\le 10\tau_{l}.
\end{equation*}
Comparing with \eqref{cpoint-equa1}-\eqref{cpoint-equa2} and using equations \eqref{multiplier-S}-\eqref{multiplier-lambda}, whenever
\begin{align*}
 \beta_{k}^{-1}\sqrt{\|S^{l,k}\!-\!\widehat{S}^{l,k-1}\|_F^2+|\lambda^{l,k}\!-\!\widehat{\lambda}^{l,k-1}|^2}\le\tau_{l},
\end{align*}
the iterate $V^{l,k}$ is a $\tau_{l}$-approximate stationary point of \eqref{gfac-penalty2} associated with $\rho=\rho_{l}$. Then, the above two inequalities can be used as the stop condition of Algorithm \ref{ALM}. 
\end{remark} 

 The following proposition shows that under a suitable condition the sequence $\{V^{l,k}\}_{k\in\mathbb{N}}$ generated by Algorithm \ref{ALM} with $\overline{\beta}=\infty$ is asymptotically feasible. Together with Remark \ref{remark-ALM} (c), we conclude that Algorithm \ref{ALM} can give arise to a $\tau_l$-approximate stationary point of \eqref{gfac-penalty2} for $\rho=\rho_{l}$ within a finite number of steps.  
 \begin{proposition}\label{ALM-feasible}
 Let $\{V^{l,k}\}_{k\in\mathbb{N}}$ be generated by Algorithm \ref{ALM} with $\overline{\beta}=\infty$. If there exists $c_{L}>0$ such that $\mathcal{L}_{\rho_l,\beta_{k-1}}(V^{l,k},\widehat{S}^{l,k-1},\widehat{\lambda}^{l,k-1})\le c_{L}$ for all $k\in\mathbb{N}$, then $\{V^{l,k}\}_{k\in\mathbb{N}}$ is bounded and its every cluster point is feasible to \eqref{gfac-penalty2} with $\rho=\rho_{l}$.
\end{proposition}
\begin{proof}
 We first prove the boundedness of $\{V^{l,k}\}_{k\in\mathbb{N}}$. From the given assumption and the expression of $\mathcal{L}_{\rho_l,\beta_{k-1}}(V^{l,k},\widehat{S}^{l,k-1},\widehat{\lambda}^{l,k-1})$, for each  $k\in\mathbb{N}$, it holds
\begin{align*}
 &\rho_{l}(\|V^{l,k}\|_F^2\!-\!\|(V^{l,k})^{\top}V^{l,k}\|_F)-\frac{1}{2\beta_{k-1}}\|\widehat{S}^{l,k-1}\|_F^2+\langle G,(V^{l,k})^{\top}V^{l,k}\rangle\\
 &+\widehat{\lambda}^{l,k-1}\big(\langle E,(V^{l,k})^{\top}V^{l,k}\rangle\!-\!1\big)+\frac{1}{2}\beta_{k-1}\big(\langle E,(V^{l,k})^{\top}V^{l,k}\rangle\!-\!1\big)^2\\
 &+\frac{1}{2}\beta_{k-1}\big\|(V^{l,k})^{\top}V^{l,k}+\beta_{k-1}^{-1}\widehat{S}^{l,k-1}-\Pi_{\mathcal{Y}}((V^{l,k})^{\top}V^{l,k}+\beta_{k-1}^{-1}\widehat{S}^{l,k-1})\big\|_F^2\le c_{L}.
 \end{align*} 
 In view of the boundedness of $\{\widehat{S}^{l,k}\}_{k\in\mathbb{N}}$, the first two terms on the left hand side of the above inequality are bounded from below. Since $\beta_k\to\infty$ as $k\to\infty$, using the above inequality and the arguments by contradiction can prove that the sequence $\{V^{l,k}\}_{k\in\mathbb{N}}$ is bounded. Then, from the definition of $Y^{l,k}$, it follows that the sequence $\{Y^{l,k}\}_{k\in\mathbb{N}}$ is bounded. Let $\overline{V}^{l}$ be an arbitrary cluster point of $\{V^{l,k}\}_{k\in\mathbb{N}}$. There exists an index set $\mathcal{K}\subset\mathbb{N}$ such that $\lim_{\mathcal{K}\ni k\to\infty}V^{l,k}=\overline{V}^{l}$. Passing the limit $\mathcal{K}\ni k\to\infty$ to the above inequality and the equality $Y^{l,k}=\Pi_{\mathcal{Y}}((V^{l,k})^{\top}V^{l,k}+\beta_{k-1}^{-1}\widehat{S}^{l,k-1})$ leads to 
\[
 \big(\langle E,(\overline{V}^{l})^{\top}\overline{V}^{l}\rangle-1\big)^2\!+\|(\overline{V}^{l})^{\top}\overline{V}^{l}-\Pi_{\mathcal{Y}}((\overline{V}^{l})^{\top}\overline{V}^{l})\|_F^2=0\ \ {\rm and}\ \ \overline{Y}^{l}=\Pi_{\mathcal{Y}}((\overline{V}^{l})^{\top}\overline{V}^{l}).
\]
This shows that $\overline{V}^{l}$ is a feasible point of \eqref{gfac-penalty2} associated with $\rho=\rho_{l}$. 
\end{proof}

Recall that our aim is to illustrate the application of global exact penalty \eqref{gfac-penalty} in designing relaxation approaches to seek a rank-one approximate feasible solution, rather than to provide a well-developed relaxation algorithm. For more theoretical analysis of Algorithm \ref{Alg-factor} and Algorithm \ref{ALM}, we leave them for the future work. 
\subsection{Numerical experiment}\label{sec5.2}  

We test the performance of Algorithm \ref{Alg-factor} armed with Algorithm \ref{ALM} (EPalm, for short) for solving the problem \eqref{QAP}. To validate its efficiency, we compare its performance with that of the commercial solver Gurobi, which uses a branch and bound method to solve the following reformulation of \eqref{QAP}:
\begin{equation*}
 \min_{x\in\{0,1\}^{p}}\left\{\langle x,Gx\rangle\ \ {\rm s.t.}\ \ 
 \begin{pmatrix}
  I\otimes e^{\top}\\  e^{\top}\otimes I  
 \end{pmatrix}x=e\right\}.
\end{equation*}
For the subsequent tests, the matrix $G$ in the penalty problem \eqref{gfac-penalty2} is scaled to $\widetilde{G}=G/\max\{1,\|G\|\}$, while the one in the above minimization keeps unscaled.   

\begin{figure}[h]
\centering
\includegraphics[width=1.0\textwidth]{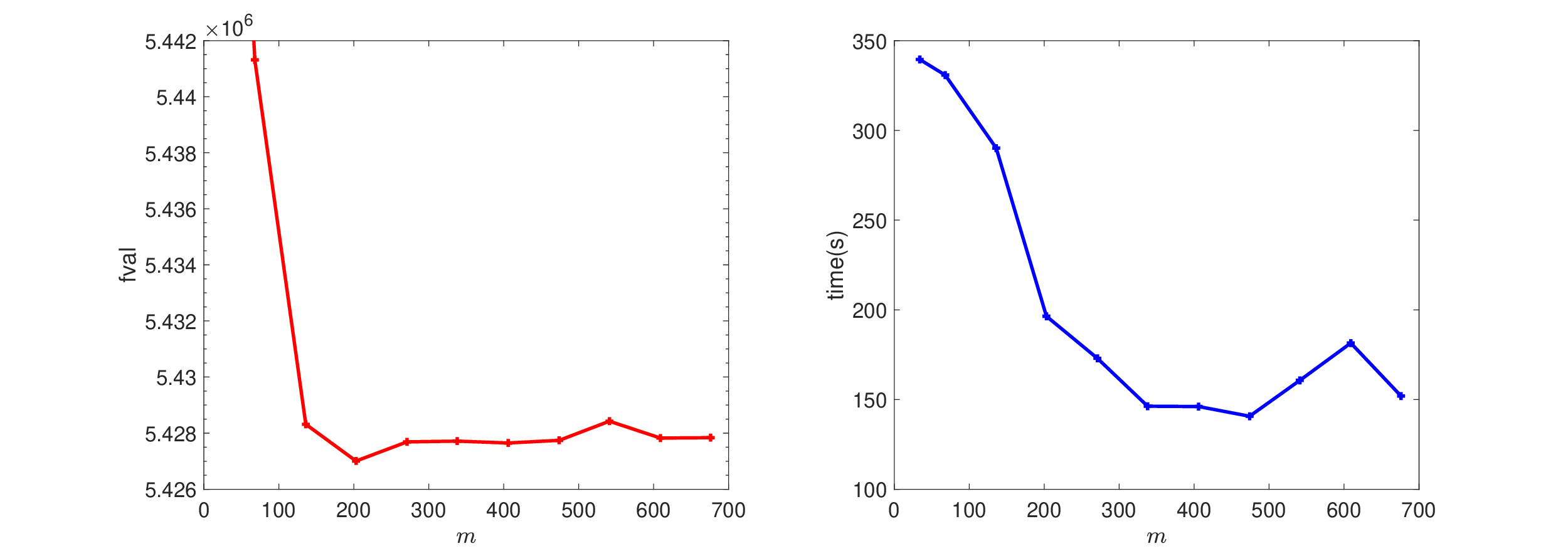}
\caption{Influence of $m$ on the objective value given by EPalm and its running time}
 \label{fig2}
\end{figure}

For the implementation of EPalm, we first take a closer look at the choice of the parameters involved in Algorithm \ref{Alg-factor}. The choice of $m$ involves a trade-off between the quality of solutions and the computation cost. Figure \ref{fig2} plots the objective value yielded by EPalm with $m\in\{\lceil \tau p\rceil\,|\,\tau=0.05,0.1,0.2,\ldots,1\}$ for solving “bur26a'' and its running time. We see that the objective values associated with $m=\lceil \tau p\rceil$ for $\tau<0.2$ is much worse than those corresponding to $m=\lceil \tau p\rceil$ for $\tau\ge 0.2$, i.e., a small $m$ (say, $m<100$) makes the penalty problem \eqref{gfac-penalty2} vulnerable to more bad stationary points; the running time decreases as $\tau$ increases from $0.05$ to $0.7$, and as $\tau$ continues to increase to $1$, it has a little increase.  For a smaller $m$, though the algorithm for solving the subproblems of ALM needs less computation cost at each iteration, it requires more iterations because the subproblems of ALM become more difficult, and the total running time becomes more. A trade-off is to choose $m=\lceil p/3\rceil$ for all tests. 

\begin{figure}[h]
\centering
\subfigure[\label{fig1a}]{\includegraphics[width=0.48\textwidth]{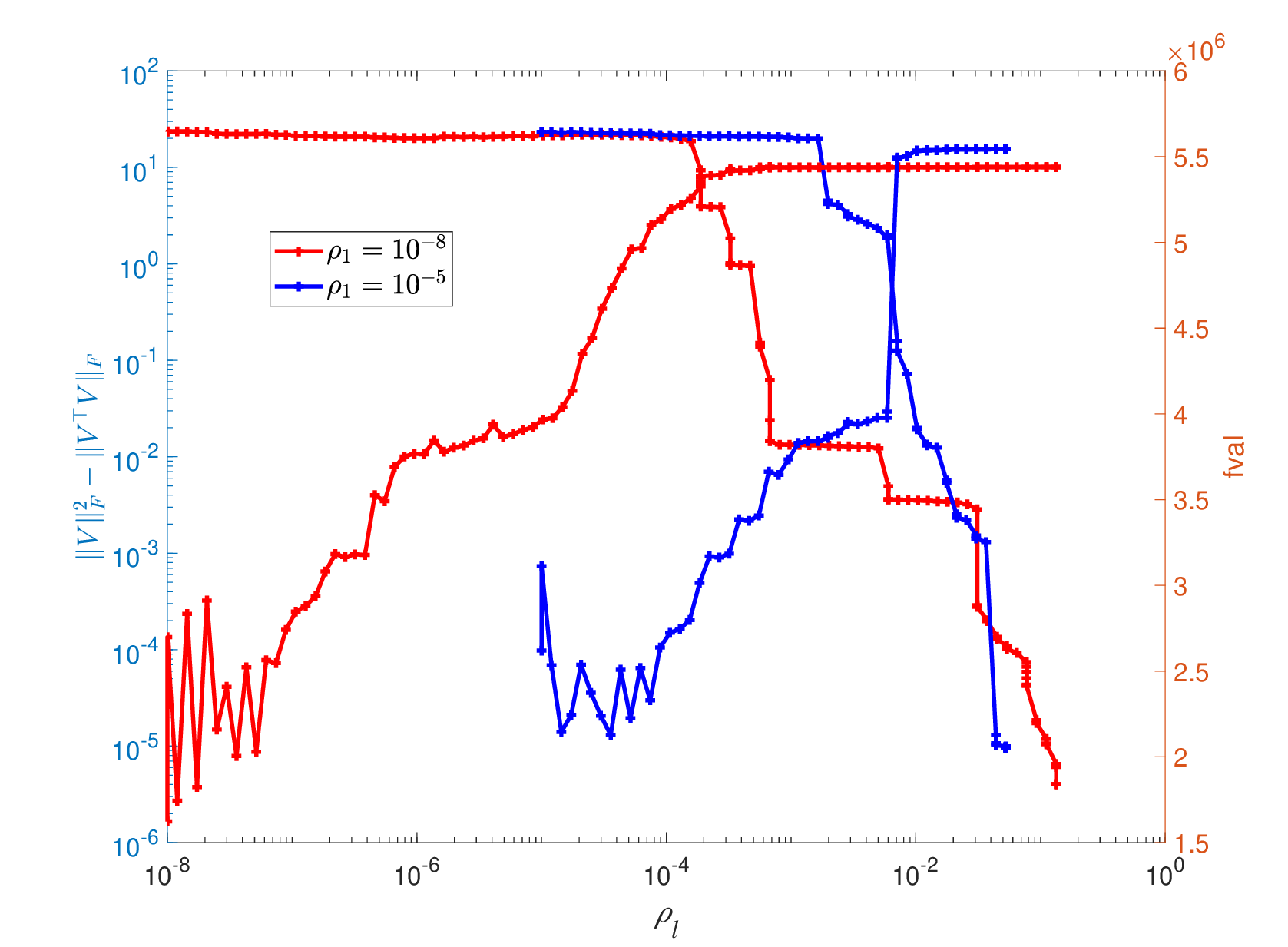}}
\subfigure[\label{fig1b}]{\includegraphics[width=0.47\textwidth]{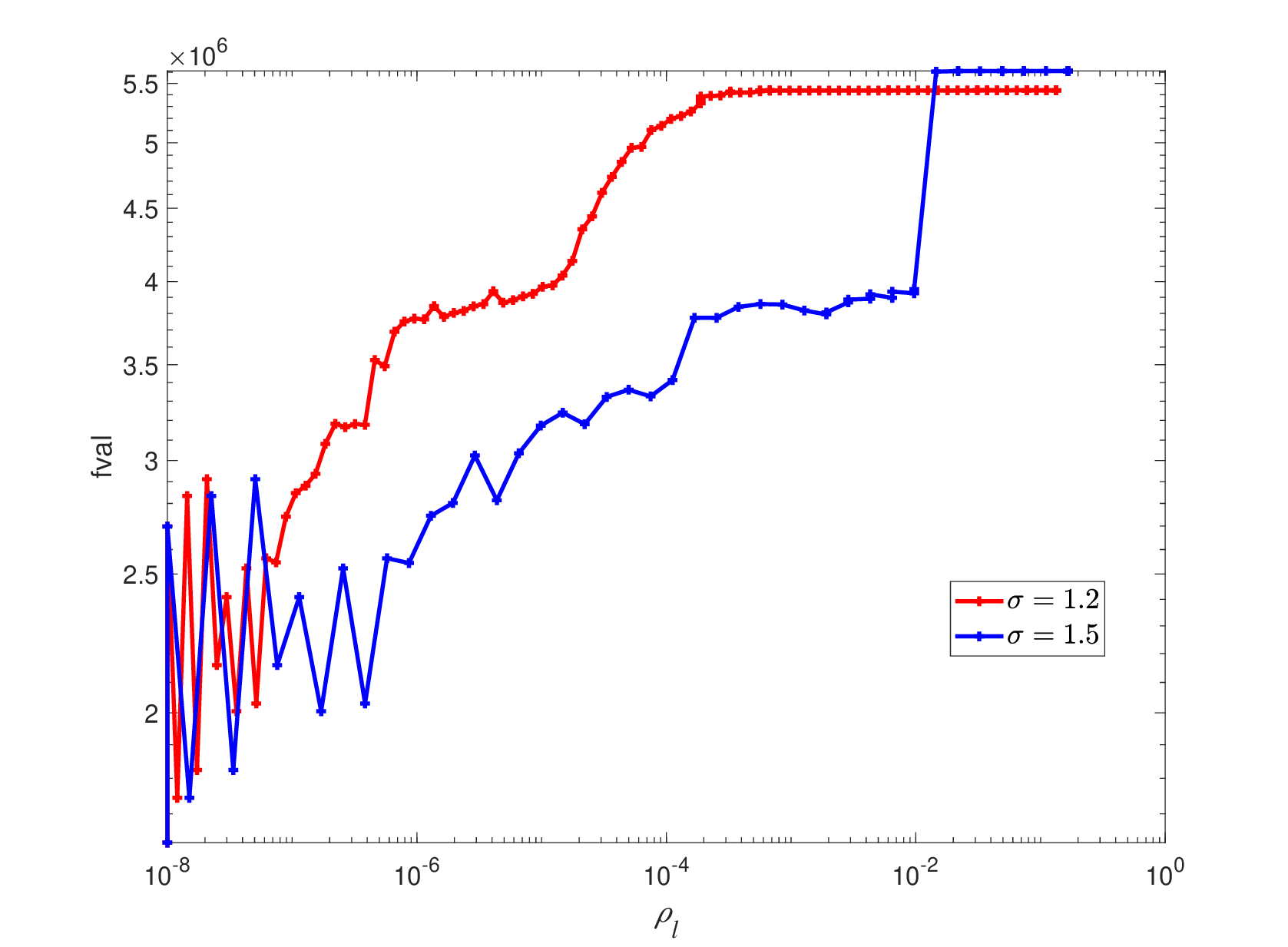}}
\caption{Influence of $\rho_1$ and $\sigma$ on the objective value and rank returned by EPalm}
\label{fig1}
\end{figure} 

Next we pay our attention to the choice of the initial penalty parameter $\rho_1$ and the adjusting coefficient $\sigma$ of the penalty parameter. Figure \ref{fig1} (a) shows that Algorithm \ref{Alg-factor} for $\sigma=1.2$ produces the lower objective value $\langle G,V^{l_f}\rangle$ and penalty term $\|V^{l_f}\|_F^2-\|(V^{l_f})^{\top}V^{l_f}\|_{F}$ with $\rho_1=10^{-8}$ than with $\rho_1=10^{-5}$. Figure \ref{fig1} (b) indicates that Algorithm \ref{Alg-factor} for $\rho_1=10^{-8}$ returns the better objective value with $\sigma=1.2$ than with $\sigma=1.5$. This shows that Algorithm \ref{Alg-factor} with the smaller $\rho_1$ and $\sigma$ leads to a better objective value, so we choose $\rho_1=10^{-8}$ and $\sigma=1.2$ for the subsequent tests. The other parameters of Algorithm \ref{Alg-factor} are chosen to be as follows
 \[
  l_{\rm max}=500,\, \rho_{\rm max}=10^5,\,\epsilon_1=\epsilon_2=10^{-5},\,\underline{\tau}=10^{-3},\,\tau_1=0.5\ \ {\rm and}\ \ \varsigma=0.9.
 \] 
 To capture a high-quality solution, it is reasonable to require the initial $V^0\in\mathbb{R}^{m\times p}$ to have a full row rank. By \cite[Corollary 5.35]{Vershynin12}, when $V^0$ is generated randomly in MATLAB command ${\rm randn}(m,p)$, there is a high probability for $V^0$ to have a full row rank. In view of this, we choose $V^0=\!\sqrt{n}V_g/\|V_g\|_F$ for the subsequent tests, where $V_g$ is generated by the MATLAB command $\textrm{randn}(m,p)$ with a seed for all test problems.

 For Algorithm \ref{ALM}, we choose $\overline{\beta}=10^5,\gamma=1.005,\beta_0=1$ and $\varepsilon=10\tau_l$. During the testing, Algorithm \ref{ALM} uses the stop condition as suggested in Remark \ref{remark-ALM} (c). We apply the limited-memory BFGS with the number of memory $l_m=15$ and the maximum number of iterations $300$ to compute the $\varepsilon$-critical point of $\mathcal{L}_{\rho_{l},\beta_k}(\cdot,\widehat{S}^{l,k},\widehat{\lambda}^{l,k})$.
 
 All the numerical tests are performed in MATLAB R2020b on a workstation running 64-bit Windows Operating System with an Intel Xeon(R) W-2245 CPU 3.90GHz and 128 GB RAM. The solution returned by EPalm is $X^*={\rm mat}(\|V^*\||Q_1^*|)$, where $V^*$ is its final iterate and $Q_1^*$ is the first column of any $Q^*\!\in\mathbb{O}^{p}((V^*)^{\top}V^*)$. The quality of $X^*$ is measured by  its relative gap and violation of feasibility, defined by 
 \[
 \textbf{gap}\!:=\frac{\rm Bval-Obj}{\rm Bval}\times 100\%
 \ \ {\rm and}\ \
 \textbf{infeas}\!:=\|(X^*)^{\top}X^*\!-I\|_F+\|\min(X^*,0)\|_F,
 \] 
 where ${\rm Bval}$ is the known best value of \eqref{QAP} and Obj is the objective value of $X^*$. We use ${\rm rgap}:=\frac{\rm Bval-rObj}{\rm Bval}\times 100\%$ to denote the relative gap of the rounded solution ${\rm round}(X^*)$, where rObj is the objective value of the rounded solution. We terminate the iteration of Gurobi once its running time exceeds that of EPalm for the same test instance or the generated objective values is not greater than ${\rm Bval}$.   
\subsubsection{Numerical results for QAPLIB data}\label{sec5.2.1}

We compare the performance of EPalm with that of Gurobi for solving \eqref{QAP} with $C=0$ and $A,B\in\mathbb{R}^{n\times n}$ from QAPLIB \cite{QAPLIB}. The \textbf{122} test examples are divided into three groups in terms of the size $n$: $n\in\!\{10,\ldots,30\},\{31,\ldots,60\}$ and $\{61,\ldots,90\}$. Tables \ref{tableQAP1}-\ref{tableQAP3} report the objective value, relative gap, and running time (in seconds) by EPalm and Gurobi, where “$-$'' in the \textbf{rObj} column means that the objective value of rounded solutions is the same as Obj. To check if EPalm yields a solution with satisfactory approximate feasibility for a suitably large $\rho$, we report the \textbf{infeas} and the final $\rho_{l_{f}}$ for EPalm. 
\setlength\tabcolsep{2.5pt}
\renewcommand\arraystretch{1.2}{
\begin{table}[h]
\setlength{\abovecaptionskip}{2pt}
\setlength{\belowcaptionskip}{0pt}
\centering
\captionsetup{font={footnotesize}}
\caption{Numerical results of EPalm and Gurobi on the small-scale QAPLIB instances}\label{tableQAP1}
\tiny
\begin{tabular}{ccccccccccccc}
\Xhline{0.55pt}\noalign{\smallskip}
{\bf No.}&{\bf Prob.}& {\bf Bval}&\multicolumn{7}{c}{\bf EPalm}
 &\multicolumn{3}{c}{\bf Gurobi}\\
			\noalign{\smallskip} \cmidrule[0.6pt](lr){4-10}\cmidrule[0.6pt](lr){11-13}
			\noalign{\smallskip}
			&   &  &Obj&rObj&gap(\%)&rgap(\%) &$\rho_{l_{\!f}}$&infeas &time(s) &Obj &gap(\%)  &time(s)\\
			\noalign{\smallskip}\Xhline{1pt}\noalign{\smallskip}
			
			1& bur26a& 5426670 & 5427506& 5427626& 1.54e-2& 1.76e-2 & 1.34e-1 & 8.1e-4 & 330.5 & 5434781 & 1.49e-1 & 330.6\\
			2& bur26b& 3817852 & 3820329& 3820275& 6.49e-2& 6.35e-2 & 1.34e-1 & 6.0e-4 & 268.4 & 3847276 & 7.71e-1 & 268.5\\
			\hline
			
			3& bur26c& 5426795 & 5428805& 5428184& 3.70e-2& 2.56e-2 & 4.48e-2 & 3.1e-4 & 149.7 & 5430053 & 6.00e-2 & 149.8\\
			4& bur26d& 3821225 & 3821629& 3821403& 1.06e-2& 4.66e-3 & 7.74e-2 & 3.3e-4 & 167.7 & 3832940 & 3.07e-1 & 167.8\\
			\hline
			
			5& bur26e& 5386879 & 5389124& 5388905& 4.17e-2& 3.76e-2 & 4.48e-2 & 4.0e-4 & 228.5 & 5402317 & 2.87e-1 & 228.5\\
			6& bur26f& 3782044 & 3782297& 3782591& 6.70e-3& 1.45e-2 & 6.45e-2 & 9.2e-4 & 322.4 & 3798840 & 4.44e-1 & 322.5\\
			\hline
			
			7 & bur26g& 10117172 & 10119241& 10119190& 2.05e-2& 1.99e-2 & 1.11e-1 & 1.3e-4 & 311.6 & 10120047 & 2.84e-2 & 311.6\\
			8 & bur26h& 7098658 & 7099195& 7099342& 7.57e-3& 9.64e-3 & 5.38e-2 & 3.5e-4 & 308.7 & 7142069 & 6.12e-1 & 308.8\\
			\hline
			9& chr12a& 9552 & 9552& --& 3.33e-3& 0 & 3.11e-2 & 4.0e-4 & 4.5 & 9552 & 0 & 0.7\\ 			 			
			10& chr12b& 9742 & 9764& 9762& 2.22e-1& 2.05e-1 & 2.59e-2 & 1.9e-4 & 5.6 & 9742 & 0 & 0.3\\
			\hline
			11&	chr12c& 11156 & 11156& --& 3.46e-3& 0 & 3.11e-2 & 5.2e-5 & 7.5 & 11156 & 0 & 0.5\\ 			
			12& chr15a& 9896 & 9896& --& 4.20e-3& 0 & 2.59e-2 & 9.9e-5 & 30.0 & 9896 & 0 & 3.5\\
			\hline
			13&	chr15b& 7990 & 7990& --& 1.08e-2& 0 & 6.45e-2 & 9.3e-5 & 9.8 & 7990 & 0 & 1.2\\ 			 			
			14 & chr15c& 9504 & 9504& --& 5.24e-3& 0 & 1.25e-2 & 3.2e-4 & 9.3 & 9504 & 0 & 0.7\\
			\hline
			15& chr18a& 11098 & 11098& --& 3.95e-3& 0 & 3.73e-2 & 6.2e-5 & 26.7 & 11098 & 0 & 9.6\\
			16& chr18b& 1534 & 1534& --& 1.32e-4& 0 & 5.38e-2 & 1.6e-4 & 12.2 & 1556 & 1.43 & 12.2\\
			\hline
			17& chr20a& 2192 & 2192& --& 7.81e-3& 0 & 2.59e-2 & 2.3e-4 & 47.1 & 2192 & 0 & 23.4\\
			18 &chr20b& 2298 & 2298& --& 6.77e-3& 0 & 3.73e-2 & 3.4e-5 & 31.1 & 2392 & 4.09 & 31.1\\
			\hline
			
			19 & chr20c& 14142 & 14151& 14150& 6.21e-2& 5.66e-2 & 6.45e-2 & 3.7e-5 & 61.6 & 14142 & 0 & 3.4\\
			20&chr22a& 6156 & 6155& 6156& -9.40e-3& 0 & 5.38e-2 & 8.9e-4 & 61.2 & 6194 & 6.17e-1 & 61.2\\
			\hline
			21& chr22b& 6194 & 6194& --& -5.08e-3& 0 & 9.29e-2 & 5.3e-4 & 82.7 & 6694 & 8.07 & 82.9\\ 			
			22& chr25a& 3796 & 3796& --& -7.62e-4& 0 & 2.59e-2 & 5.0e-4 & 115.9 & 4558 & 20.07 & 115.9\\
			\hline
			
			23& els19 & 17212548 & 17215782& 17212548& 1.88e-2& 0 & 4.48e-2 & 2.9e-4 & 28.8 & 17212548 & 0 & 7.3\\
			24&esc16a& 68 & 68& --& 3.56e-4& 0 & 7.74e-2 & 1.7e-5 & 31.0 & 68 & 0 & 0.4\\
			\hline
			
			25&	esc16b& 292 & 292& --& 3.20e-3& 0 & 7.74e-2 & 7.6e-5 & 40.0 & 292 & 0 & 0.6\\
			26&	esc16c& 160 & 160& --& 2.18e-3& 0 & 2.59e-2 & 6.9e-5 & 36.8 & 160 & 0 & 1.6\\
			\hline
			
			27 & esc16d& 16 & 16& --& -1.88e-2& 0 & 5.38e-2 & 8.7e-4 & 41.7 & 16 & 0 & 0.4\\
			28&	esc16e& 28 & 28& --& -9.89e-3& 0 & 4.48e-2 & 5.4e-5 & 10.4 & 28 & 0 & 0.2\\
			\hline

            29 & esc16f & 0 & 0& --& 0& 0 & 2.59e-2 & 2.8e-4 & 4.2 & 0 & 0 & 0.02\\
			30&	esc16g& 26 & 26& --& 7.81e-4& 0 & 7.74e-2 & 1.1e-4 & 13.7 & 26 & 0 & 0.2\\
			31&	esc16h& 996 & 996& --& 8.90e-4& 0 & 5.38e-2 & 1.5e-4 & 12.9 & 996 & 0 & 1.1\\
			\hline
			
			32&esc16i& 14 & 14& --& 2.36e-3& 0 & 4.48e-2 & 6.8e-5 & 48.3 & 14 & 0 & 0.1\\
			33&esc16j& 8 & 8& --& -2.97e-4& 0 & 4.48e-2 & 2.2e-4 & 17.7 & 8 & 0 & 0.1\\
			\hline
			
			34 & had12& 1652 & 1652& --& 1.88e-3& 0 & 4.48e-2 & 9.8e-5 & 9.2 & 1666 & 8.47e-1 & 9.2\\
			35 & had14& 2724 & 2724& --& 4.48e-4& 0 & 6.45e-2 & 6.7e-5 & 17.0 & 2724 & 0 & 4.2\\
			\hline
			
			36 & had16& 3720 & 3720& --& 8.62e-3& 0 & 7.74e-2 & 1.2e-4 & 23.6 & 3728 & 2.15e-1 & 23.7\\
			37 & had18& 5358 & 5358& --& -3.05e-3& 0 & 5.38e-2 & 4.5e-4 & 44.1 & 5398 & 7.47e-1 & 44.2\\	
			38 & had20& 6922 & 6922& --& -2.35e-4& 0 & 9.29e-2 & 2.7e-4 & 63.4 & 6964 & 6.07e-1 & 63.5\\
			\hline
			
			39& kra30a& 88900 & 88905& 88900& 5.11e-3& 0 & 5.38e-2 & 3.4e-4 & 329.5 & 91300 & 2.70 & 329.5\\
			40&	kra30b& 91420 & 93024& 93030& 1.75& 1.76 & 7.74e-2 & 4.6e-4 & 326.3 & 93960 & 2.78 & 326.4\\
			\hline
			41&	lipa20a& 3683 & 3683& --& 2.18e-4& 0 & 3.73e-2 & 2.5e-4 & 18.4 & 3797 & 3.10 & 18.5\\
			42&	lipa20b& 27076 & 27078& 27076& 5.93e-3& 0 & 1.25e-2 & 3.2e-4 & 11.0 & 27076 & 0 & 0.9\\
			\hline
			43&lipa30a& 13178 & 13177& 13178& -7.29e-3& 0 & 4.48e-2 & 9.1e-4 & 150.1 & 13438 & 1.97 & 150.1\\
			44&lipa30b& 151426 & 151428& 151426& 1.08e-3& 0 & 7.23e-3 & 2.8e-4 & 55.8 & 151426 & 0 & 5.1\\
			
			\hline
			\Xhline{0.55pt}
	\end{tabular}
	\end{table}
}
\setcounter{table}{0}
\setlength\tabcolsep{2.4pt}
\renewcommand\arraystretch{1.2}{
  \begin{table}[h]
		\setlength{\abovecaptionskip}{2pt}
		\setlength{\belowcaptionskip}{0pt}
		\centering
		\captionsetup{font={footnotesize}}
		\caption{(continued)}\label{tableQAP1c}
		\tiny
		\begin{tabular}{ccccccccccccc}
			\Xhline{0.55pt}\noalign{\smallskip}
			{\bf No.}&{\bf Prob.}& {\bf Bval}&\multicolumn{7}{c}{\bf EPalm}
			&\multicolumn{3}{c}{\bf Gurobi}\\
			\noalign{\smallskip} \cmidrule[0.6pt](lr){4-10}\cmidrule[0.6pt](lr){11-13}\noalign{\smallskip}
			&  &  & Obj&rObj&gap(\%)&rgap(\%)&$\rho_{l_{\!f}}$&infeas &time(s) &Obj &gap(\%) &time(s)\\
			\noalign{\smallskip}\Xhline{1pt}\noalign{\smallskip}
			
			45& nug12& 578 & 578& --& -7.25e-3 & 0 & 5.38e-2 & 4.5e-4 & 93.3 & 578 & 0 & 16.8\\
			46& nug14& 1014 & 1018& --& 3.95e-1 & 3.94e-1 & 4.48e-2 & 3.2e-5 & 25.4 & 1030 & 1.58 & 25.4\\	
			47&	nug15& 1150 & 1152& --& 1.72e-1 & 1.74e-1 & 5.38e-2 & 1.3e-4 & 25.7 & 1180 & 2.61 & 25.7\\
			\hline

			48& nug16a& 1610 & 1610& --& -5.69e-3 & 0 & 6.45e-2 & 5.0e-4 & 53.0 & 1664 & 3.35 & 53.0\\
			49&	nug16b& 1240 & 1258& --& 1.45 & 1.45 & 7.74e-2 & 3.5e-4 & 34.7 & 1294 & 4.35 & 34.7\\
			\hline
			
			50&	nug17& 1732 & 1734& --& 1.16e-1 & 1.15e-1 & 1.50e-2 & 1.3e-4 & 43.6 & 1758 & 1.50 & 43.6\\
			51& nug18& 1930 & 1944& --& 7.26e-1 & 7.25e-1 & 3.11e-2 & 9.6e-5 & 47.0 & 1976 & 2.38 & 47.0\\
			\hline
			
			52& nug20& 2570 & 2596& --& 1.02 & 1.01 & 2.59e-2 & 1.7e-4 & 66.2 & 2662 & 3.58 & 66.3\\	
			53& nug21& 2438 & 2452& --& 5.80e-1 & 5.74e-1 & 3.73e-2 & 2.4e-4 & 99.9 & 2508 & 2.87 & 99.9\\
			\hline
			
			54&	nug22& 3596 & 3658& --& 1.72 & 1.72 & 6.45e-2 & 1.9e-5 & 141.4 & 3664 & 1.89 & 141.5\\
			55&	nug24& 3488 & 3512& --& 6.91e-1 & 6.88e-1 & 3.73e-2 & 3.1e-4 & 148.8 & 3600 & 3.21 & 149.0\\
			\hline
			
			56&	nug25& 3744 & 3750& --& 1.64e-1 & 1.60e-1 & 3.73e-2 & 1.9e-4 & 170.9 & 3810 & 1.76 & 171.0\\
			57& nug27& 5234 & 5238& --& 6.90e-2 & 7.64e-2 & 4.48e-2 & 1.0e-3 & 241.4 & 5430 & 3.74 & 241.5\\
			
			\hline

			58&	nug28& 5166 & 5248& --& 1.59 & 1.59 & 7.74e-2 & 3.0e-4 & 270.8 & 5510 & 6.66 & 270.9\\
			59& nug30& 6124 & 6182& --& 9.46e-1 & 9.47e-1 & 3.11e-2 & 2.1e-4 & 337.2 & 6380 & 4.18 & 337.3\\
			\hline

			60& rou12& 235528 & 235534& 235528& 2.62e-3 & 0 & 3.73e-2 & 4.7e-5 & 56.1 & 238134 & 1.11 & 56.1\\
			61& rou15& 354210 & 354197& 354210& -3.62e-3 & 0 & 9.29e-2 & 2.2e-4 & 103.4 & 360278 & 1.71 & 103.4\\
			62& rou20& 725522 & 731059& 730970& 7.63e-1 & 7.51e-1 & 4.48e-2 & 4.3e-4 & 326.2 & 732182 & 9.18e-1 & 326.2\\
			\hline
			
			63& scr12& 31410 & 31416& 31410& 1.84e-2 & 0 & 3.11e-2 & 2.9e-4 & 7.0 & 31410 & 0 & 4.6\\
			64& scr15& 51140 & 51139& 51140& -2.02e-3 & 0 & 6.45e-2 & 4.6e-4 & 8.8 & 51140 & 0 & 8.8\\
			65&	scr20& 110030 & 110032& 110030& 1.80e-3 & 0 & 2.59e-2 & 4.2e-4 & 85.6 & 110994 & 8.76e-1 & 85.6\\
			\hline
			
			66&	tai10a& 135028 & 135013& 135028& -1.08e-2 & 0 & 2.16e-2 & 6.3e-4 & 13.8 & 135028 & 0 & 0.2\\	
			67& tai12a& 224416 & 224420& 224416& 1.59e-3 & 0 & 7.74e-2 & 3.6e-5 & 3.4 & 236920 & 5.57 & 3.4\\
			\hline
			
			68& tai12b& 39464925 & 39477370& 39477247& 3.15e-2 & 3.12e-2 & 1.50e-2 & 2.1e-4 & 20.2 & 39464925 & 0 & 20.1\\
			69& tai15a& 388214 & 388872& 388870& 1.69e-1 & 1.69e-1 & 6.45e-2 & 1.3e-5 & 131.7 & 397376 & 2.36 & 131.8\\
			\hline
			
			70&	tai15b& 51765268 & 52254153& 52251783& 9.44e-1 & 9.40e-1 & 1.50e-2 & 2.4e-4 & 10.4 & 51825455 & 1.16e-1 & 10.5\\
			71& tai17a& 491812 & 493659& 493662& 3.76e-1 & 3.76e-1 & 3.11e-2 & 2.2e-4 & 161.0 & 510114 & 3.72 & 161.0\\	
			\hline
			
			72&	tai20a& 703482 & 710788& 710786& 1.04 & 1.04 & 6.45e-2 & 7.2e-5 & 400.5 & 734718 & 4.44 & 400.5\\			
			73& tai20b& 122455319 & 122453700& 122455319& -1.32e-3 & 0 & 6.45e-2 & 4.2e-4 & 127.9 & 123143224 & 5.62e-1 & 127.9\\
			\hline
			
			74& tai25a& 1167256 & 1191771& 1191798& 2.10 & 2.10 & 2.59e-2 & 4.8e-4 & 289.1 & 1218210 & 4.37 & 289.2\\
			75& tai25b& 344355646 & 348335422& 348328772& 1.16 & 1.15 & 5.38e-2 & 2.0e-4 & 438.9 & 349538768 & 1.51 & 438.9\\
			\hline
			
			76& tai30a& 1818146 & 1863382& 1863600& 2.49 & 2.50 & 2.59e-2 & 1.1e-3 & 287.1 & 1907118 & 4.89 & 287.1\\ 			 			
			77& tai30b& 637117113 & 644179202& 644176415& 1.11 & 1.11 & 9.29e-2 & 2.9e-4 & 965.8 & 650340150 & 2.08 & 965.9\\
			\hline
			
			78& tho30& 149936 & 151476& --& 1.03 & 1.03 & 4.48e-2 & 1.7e-4 & 350.8 & 154014 & 2.72 & 350.9\\
			\hline
			\Xhline{0.55pt}
		\end{tabular}
	\end{table}
}

\begin{figure}[h]
\centering
\includegraphics[width=1\textwidth]{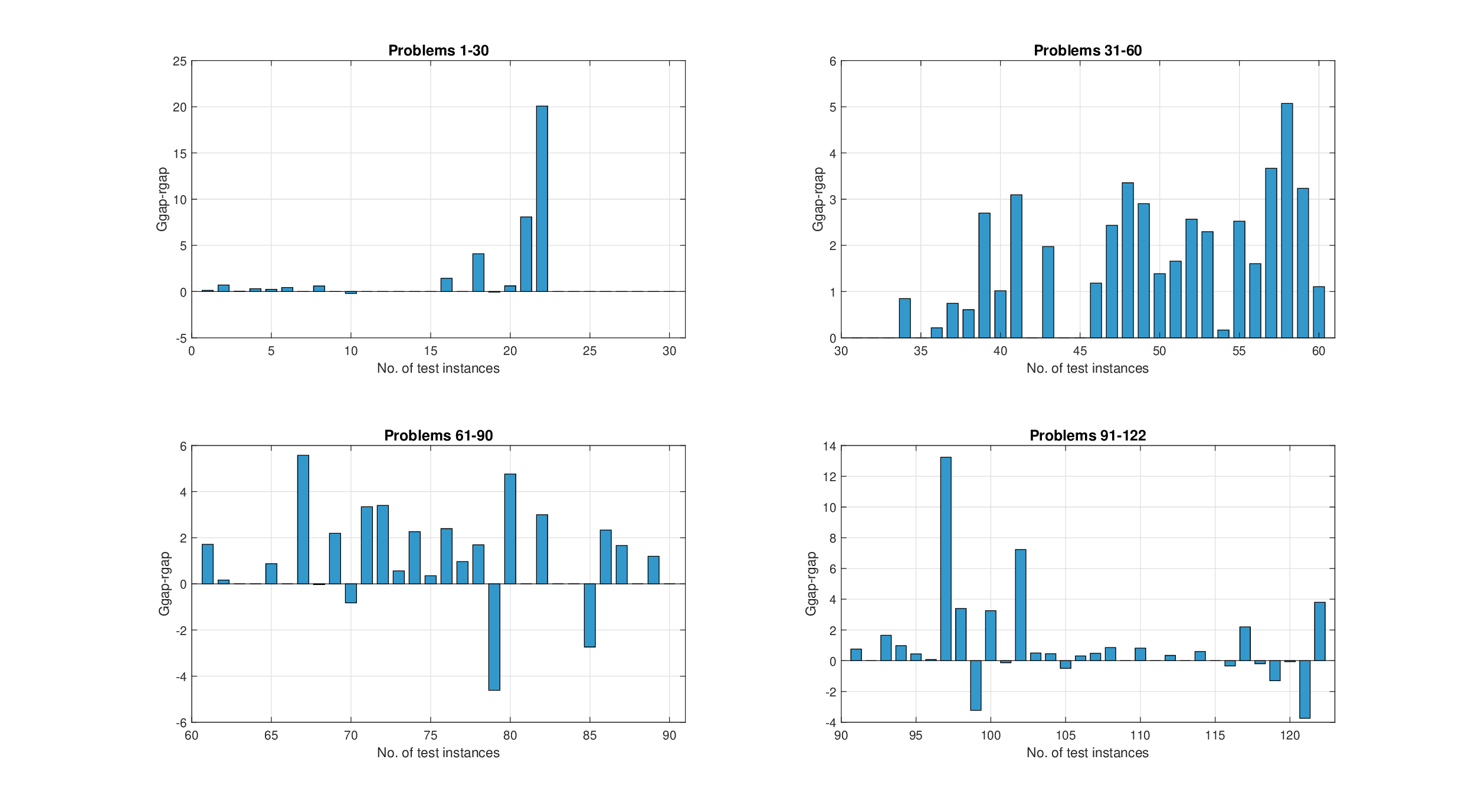}
\caption{The gap difference of Gurobi from that of EPalm for the QAPLIB instances (where Ggap appearing in ordinate is the one yielded by Gurobi)}
 \label{fig3}
\end{figure}

From Tables \ref{tableQAP1}-\ref{tableQAP3}, the violation of feasibility yielded by EPalm is about $10^{-4}$ for all instances, which has no influence on the rounding of solutions though the values in the \textbf{rgap} (resp. \textbf{rObj}) column are different from those in the \textbf{gap} (resp. \textbf{Obj}) column. For every test instance, the final penalty factor $\rho$ is not greater than $\rho_{l_{f}}\|G\|$, which is consistent with the global exact penalty results in Section \ref{sec4} by virtue of the values in the $\rho_{l_{f}}$ column. From Tables \ref{tableQAP1}-\ref{tableQAP2}, among the $\textbf{108}$ small and medium-scale examples, there are $\textbf{68}$ and $\textbf{9}$ ones, respectively, for which EPalm yields lower and higher objective values than Gurobi; the maximum rgap by EPalm is $\textbf{10.77}\%$, while the one by Gurobi is $\textbf{20.07}\%$; there are $\textbf{8}$ examples by EPalm and  $\textbf{36}$ examples by Gurobi, respectively, whose rgap is more than $\textbf{2.5}\%$. For the $\textbf{14}$ large-scale instances in Table \ref{tableQAP3}, there are $\textbf{5}$ one for which EPalm produces better objective values and worse objective values, respectively, and the maximum rgap $\textbf{8.24}\%$ by EPalm is greater than that of Gurobi $\textbf{5.93}\%$. To sum up, EPalm is superior to Gurobi in terms of the relative gap and the number of examples with better objective values. The gap difference of Gurobi from that of EPalm in Figure \ref{fig3} also demonstrates this. Observe that the running time of Gurobi for those instances with zero gap is less than that of EPalm, since the former uses the known best values as either of the stop conditions.           
 
\setlength\tabcolsep{2.3pt}
\renewcommand\arraystretch{1.2}{
	\begin{table}[htbp]
		\setlength{\abovecaptionskip}{2pt}
		\setlength{\belowcaptionskip}{0pt}
		\centering
		\captionsetup{font={footnotesize}}
		\caption{Numerical results of EPalm and Gurobi on the medium-scale QAPLIB instances}\label{tableQAP2}
		\tiny
		\begin{tabular}{cccccccccccccc}
			\Xhline{0.55pt}\noalign{\smallskip}
			{\bf No.}&{\bf Prob.}& {\bf Bval}&\multicolumn{7}{c}{\bf EPalm}
			&\multicolumn{3}{c}{\bf Gurobi}\\
			\noalign{\smallskip} \cmidrule[0.6pt](lr){4-10}\cmidrule[0.6pt](lr){11-13}\noalign{\smallskip}
			&  &  & Obj& rObj&gap(\%)&rgap(\%)&$\rho_{l_{\!f}}$&infeas &time(s) &Obj &gap(\%) &time(s)\\
            \noalign{\smallskip}\Xhline{1pt}\noalign{\smallskip}
			
			79&esc32a& 130 & 144& --& 10.77& 10.77& 7.74e-2&1.5e-4 & 1199.5 & 138& 6.15& 1199.6\\
			80&	esc32b& 168 & 184& --& 9.52& 9.52& 2.59e-2&3.6e-4 & 453.9 & 192& 14.29& 453.9\\
			\hline
			
			81&	esc32c& 642 & 642& --&-4.13e-3& 0& 7.74e-2&3.5e-4 & 824.7 & 642& 0& 9.6\\
			82&	esc32d& 200 & 200& --& -3.38e-4& 0& 3.73e-2&2.9e-4 & 820.4 & 206& 3.00& 820.4\\
			\hline
			
			83&	esc32e& 2 & 2& --& 2.16e-2& 0& 6.45e-2&1.8e-4 & 340.9 & 2& 0& 0.2\\
			84&	esc32g& 6 & 6& --& 1.92e-2& 0& 1.18e-2&1.2e-4 & 162.4 & 6& 0& 0.2\\
			\hline
			
			85&	esc32h& 438 & 452& --&3.20& 3.20& 1.61e-1&3.9e-4 & 295.6 & 440& 4.57e-1& 295.6\\
			86& kra32& 88700 & 89401& 89400& 7.91e-1& 7.89e-1& 6.45e-2&1.5e-4 & 445.3 & 91470& 3.12& 445.3\\
			\hline
			
			87& lipa40a& 31538 & 31537& 31538&-2.43e-3& 0& 2.16e-2&5.8e-4 & 814.9 & 32063& 1.66& 815.0\\
			
			88&	lipa40b& 476581 & 476603& 476581& 4.61e-3& 0& 5.38e-2&2.5e-4 & 227.3 & 476581& 0& 18.8\\
			\hline
			
			89&lipa50a& 62093 & 62094& 62093& 2.15e-3& 0& 4.48e-2&2.8e-4 & 2866.8 & 62836& 1.20& 2867.2\\			
			90& lipa50b& 1210244 & 1210138& 1210244& -8.78e-3& 0& 4.48e-2&1.4e-3 & 748.6 & 1210244& 0& 749.0\\
			\hline
			91& lipa60a& 107218 & 107599& --& 3.55e-1& 3.55e-1& 3.73e-2&3.8e-4 & 6881.7 & 108417& 1.12& 6882.5\\						
			92& lipa60b& 2520135 & 2519883& 2520135& -9.99e-3& 0& 8.68e-3&1.1e-3 & 2291.2 & 2520135& 0& 2291.6\\
			\hline
			
			93&   sko42& 15812 & 16016& 16018& 1.29& 1.30& 3.73e-2&1.3e-3 & 1255.6 & 16280& 2.96& 1255.9\\						
			94&	 sko49& 23386 & 23813& 23812& 1.83& 1.82& 1.11e-1&3.7e-4 & 3703.5 & 24042& 2.81& 3703.8\\
			95&  sko56& 34458 & 34993& 35002& 1.55& 1.58& 9.29e-2&3.1e-3 & 6435.8 & 35154& 2.02& 6436.2\\
			\hline
			
			96& ste36a& 9526 & 9612& --& 9.01e-1& 9.03e-1& 2.16e-2&4.2e-4 & 804.4 & 9620& 9.87e-1& 804.4\\
			97&	ste36b& 15852 & 15855& 15852& 1.93e-2& 0& 4.48e-2&2.8e-4 & 576.9 & 17950& 13.23& 577.0\\
			98&	ste36c& 8239110 & 8239596& 8239110& 5.90e-3& 0& 3.11e-2&5.0e-4 & 757.1 & 8518962& 3.40& 757.3\\
			\hline
			
			99& tai35a& 2422002 & 2598272& 2598274& 7.28& 7.28& 3.33e-1&3.9e-5 & 1296.2 & 2520188& 4.05& 1296.3\\
			100& tai35b& 283315445 & 284884794& 284889077& 5.54e-1& 5.55e-1& 9.29e-2&2.6e-4 & 1789.7 & 294113088& 3.81& 1789.8\\
			\hline
			101& tai40a& 3139370 & 3296048& 3295814& 4.99& 4.98& 4.48e-2&4.7e-4 & 570.4 & 3291844& 4.86& 570.6\\					
			102& tai40b& 637250948 & 639705468& 639707985& 3.85e-1& 3.86e-1& 1.93e-1&3.1e-4 & 2968.8 & 685786663& 7.62& 2968.9\\
			\hline		
			
			103& tai50a& 4941410 & 5208364& 5208378& 5.40& 5.40& 3.73e-2&4.5e-4 & 1843.4 & 5233298& 5.91& 1843.7\\
			
			104& tai50b& 458821517 & 467695162& 467693852& 1.93& 1.93& 2.31e-1&2.0e-4 & 8656.2 & 469784387& 2.39& 8645.4\\
			\hline			
			
			105& tai60a& 7208572 & 7634772& 7634834& 5.91& 5.91& 1.50e-2&3.6e-4 & 4689.3 & 7599342& 5.42& 4689.6\\
			106& tai60b& 608215054 & 614884764& 614880454& 1.10& 1.10& 2.77e-1&1.9e-4 & 19343.7 & 616764142& 1.41& 19344.2\\
			\hline
			107& tho40& 240516 & 245471& 245532& 2.06& 2.09& 3.11e-2&3.0e-3 & 1257.6 & 246698& 2.57& 1257.7\\
			108& wil50& 48816 & 49147& 49148& 6.79e-1& 6.80e-1& 1.11e-1&2.8e-4 & 2721.5 & 49568& 1.54& 2721.9\\

			\Xhline{0.55pt}
		\end{tabular}
	\end{table}
}
\setlength\tabcolsep{2.25pt}
 \renewcommand\arraystretch{1.2}{
  \begin{table}[htbp]
		\setlength{\abovecaptionskip}{2pt}
		\setlength{\belowcaptionskip}{0pt}
		\centering
		\captionsetup{font={footnotesize}}
		\caption{Numerical results of EPalm and Gurobi on the large-scale QAPLIB instances}\label{tableQAP3}
		\tiny
		\begin{tabular}{ccccccccccccc}
		 \Xhline{0.5pt}\noalign{\smallskip}
		 {\bf No.}&{\bf Prob.}& {\bf Bval}&\multicolumn{7}{c}{\bf EPalm}
			&\multicolumn{3}{c}{\bf Gurobi}\\
			\noalign{\smallskip} \cmidrule[0.5pt](lr){4-10}\cmidrule[0.5pt](lr){11-13}\noalign{\smallskip}
			& &   & Obj& rObj&gap(\%)&rgap(\%)&$\rho_{l_{\!f}}$&infeas &time(s) &Obj &gap(\%) &time(s)\\
			\noalign{\smallskip}\Xhline{1pt}\noalign{\smallskip}

			109& esc64a& 116 & 116& --& 1.13e-3&  0&   4.48e-2& 8.61e-4& 46946.3& 116& 0& 9.0\\
			\hline
			110& lipa70a& 169755 & 170081& 170080& 1.92e-1& 1.91e-1& 3.73e-2& 2.20e-4& 19508.6& 171476& 1.01& 19508.8\\
			111& lipa70b& 4603200 & 4603229& 4603200& 6.25e-4& 0& 8.68e-3& 4.03e-4& 4559.8& 4603200& 0& 4560.1\\
			\hline
			112& lipa80a& 253195 & 254736& 254732& 6.09e-1& 6.07e-1& 3.37e-2& 4.14e-4& 18351.7& 255604& 9.51e-1& 18352.3\\
			113& lipa80b& 7763962 & 7763949& 7763962& -1.73e-4& 0& 1.25e-2& 4.05e-4& 10962.3& 7763962& 0& 705.0\\
			\hline
			114&lipa90a& 360630 & 361643& 361640& 2.81e-1& 2.80e-1& 9.29e-2& 3.22e-4& 59852.8& 363802& 8.80e-1& 59856.6\\
			115& lipa90b& 12490441 & 12490539& 12490441& 7.82e-4&  0&   6.45e-2& 4.90e-4& 18845.4& 12490441& 0& 1349.6\\
			\hline
			116& sko64& 48498 & 49539& 49546&  2.15& 2.16&  4.48e-2& 2.19e-3& 11390.2& 49382& 1.82& 11390.9\\
			117& sko72& 66256 &  67310& --& 1.59& 1.59& 2.59e-2& 1.51e-4& 17944.9& 68770& 3.79& 17945.4\\
			
			\hline
			118& sko81& 90998 &  92805& 92804& 1.99& 1.98& 1.11e-1& 6.43e-5& 34642.5& 92628& 1.79& 34643.2\\				
			119& sko90& 115534 &  119837& 119838& 3.72& 3.73& 2.59e-2& 2.77e-4& 54336.7& 118332& 2.42& 54338.8\\
			\hline
			
			120& tai64c& 1855928 & 1860728& 1860702& 2.59e-1& 2.57e-1& 7.74e-2& 2.89e-4& 28081.1& 1859480& 1.91e-1& 28081.7\\
			
			121& tai80a& 13557864 & 14675229& 14675402& 8.24& 8.24 & 3.73e-2 & 3.24e-4& 19982.6& 14167534& 4.50& 19983.1\\
			122& tai80b& 818415043 & 835807445& 835808411& 2.13& 2.13& 1.93e-1 & 1.63e-4& 50554.8& 866948331& 5.93& 50556.3\\
			
		\Xhline{0.50pt}
		\end{tabular}
		\end{table}
       }

We also compare the objective values obtained by imposing the rounding on the solutions by EPalm with the upper bounds yielded by rPRSM of \cite{Graham22}. For the $\textbf{84}$ test instances listed in \cite[Tables 1-3]{Graham22}, EPalm yields the better upper bounds than rPRSM for $\textbf{53}$ instances, and the worse upper bounds only for $\textbf{3}$ instances.

\subsubsection{Numerical results for Dre problems}\label{sec5.2.2}

As described in \cite{Drezner05}, the `dre' test instances are based on a rectangular grid where all nonadjacent nodes have zero weight. This way, a pair exchange of the optimal permutation will result in many adjacent pairs becoming non-adjacent, which making the objective value will increase quite steeply. The `dre' instances are difficult to solve, especially for many metaheuristic-based methods because they 
are ill-conditioned and hard to break out the basin of the local minimal. These instances are available in \url{http://business.fullerton.edu/zdrezner}, and the known best solutions for them have been found by branch and bound in \cite{Drezner05}. Table \ref{tabledre} reports the results of EPalm and Gurobi for the $\textbf{10}$ `dre' test instances. EPalm produces the best objective values for all instances, but Gurobi returns the best ones only for $\textbf{4}$ instances and the relative gaps for the other $\textbf{6}$ ones are more than $50\%$.  

\setlength\tabcolsep{4pt}
\renewcommand\arraystretch{1.2}{
 \begin{table}[h]
 \setlength{\abovecaptionskip}{2pt}
 \setlength{\belowcaptionskip}{0pt}		
 \centering
 \captionsetup{font={footnotesize}}
 \caption{Numerical results of EPalm and Gurobi on the 10 `dre' instances}\label{tabledre}
 \tiny
 \begin{tabular}{ccccccccccccc}
 \Xhline{0.55pt}\noalign{\smallskip}
 {\bf No.} &{\bf Prob. }& {\bf Bval} &\multicolumn{7}{c}{\bf EPalm}
 &\multicolumn{3}{c}{\bf Gurobi}\\
 \noalign{\smallskip}\cmidrule[0.6pt](lr){4-10}\cmidrule[0.6pt](lr){11-13}\noalign{\smallskip}
 &  &  & Obj& rObj &gap(\%)&rgap(\%)&$\rho_{l_{\!f}}$&infeas &time(s) &Obj&gap(\%) &time(s) \\
 \noalign{\smallskip}\Xhline{1pt}\noalign{\smallskip}
 			1 & dre15& 306 & 306& --& 5.32e-2 & 0& 1.50e-2& 3.5e-4 & 26.4 & 306& 0& 7.0\\
 			2 & dre18& 332 & 332& --& 2.66e-3 & 0& 2.59e-2& 1.9e-4 & 37.7 & 332& 0& 11.1\\
 			\hline
 			
 			3& dre21& 356 & 356& --& 6.10e-3 & 0& 1.80e-2& 2.4e-4 & 49.3 & 356& 0& 31.8\\
 			4& dre24& 396 & 396& --& 1.83e-2 & 0& 6.03e-3& 5.6e-4 & 95.9 & 396& 0& 82.1\\
 			\hline
 			
 			5& dre28& 476 & 476& --& 9.18e-4 & 0& 4.48e-2& 7.0e-5 & 134.0 & 782& 64.29& 134.1\\
 			6& dre30& 508 & 508& --& 2.92e-3 & 0& 3.11e-2& 1.4e-4 & 159.3 & 770& 51.57& 159.3\\
 			\hline
 			
 			7&dre42& 764 & 764& --& 5.35e-3 & 0& 3.73e-2& 2.9e-4 & 807.1 & 1344& 75.92& 807.2\\
 			8 &dre56& 1086 & 1086& --& -6.09e-4 & 0& 3.11e-2& 2.1e-4 & 3966.7 & 2246& 106.81& 3966.8\\
 			\hline
 			
 			9& dre72& 1452 & 1452& --& 1.50e-3 & 0& 1.80e-2& 1.6e-4 & 20359.6 & 2968& 104.41& 20360.1\\
 			10&	dre90& 1838 & 1838& --& -1.08e-3& 0& 3.73e-2& 2.9e-4 &52167.6 & 3922& 113.84& 52169.5\\
 			\hline 			
 			\Xhline{0.55pt}
 		\end{tabular}
 	\end{table}
 }
\section{Conclusion}\label{sec6}

 We presented three equivalent rank-one DNN reformulations for the QAP, including the one considered in \cite{Jiang21}, and established the locally Lipschitzian error bounds for their feasible set $\Gamma$. With the help of these error bounds, the penalty problem \eqref{gDNN-penalty} induced by the DC reformulation of the rank-one constraint was proved to be a global exact penalty of \eqref{gDNN-mpec}, and the same conclusion was achieved for their BM factorizations \eqref{gfac-penalty} and \eqref{gfac-mpec}. This not only recovers the global exact penalty result in \cite{Jiang21} without the calmness assumption on $\Upsilon_0$, but also enriches greatly the global exact penalty results for the rank-one DNN reformulations of the QAP. We also proposed a relaxation approach with the exact penalty \eqref{gfac-penalty2} to illustrate their application in designing continuous relaxation methods for seeking rank-one approximate feasible solutions. This approach yields a rank-one approximate feasible solution by searching for a finite number of approximate stationary points of the penalty subproblems with increasing $\rho$ via the ALM. Numerical comparisons with the commercial solver Gurobi showed that the proposed EPalm has an advantage over Gurobi in terms of the quality of solutions, though its running time for the large-scale instances are still too long. We leave this challenging topic for our future work.

\bibliographystyle{siamplain}
\bibliography{references}            
\end{document}